**Topics In Elliptic Curves Cryptography: The Groups of Points**
**N.A. Carella, 2011**

# Content

# 1. Introduction

This note is concerned with the elementary theory and practice of elliptic curves cryptography, the new generation of public key systems. The material and coverage are focused on the groups of points of elliptic curves and algebraic curves, but not exclusively. The level of detail ranges from the very simple to the very difficult and open problems. The detailed proofs of many simple results are included or discussed. However, the long and difficult proofs of difficult results are not included, but relevant references are given. The reader should consult the literature for further elaboration on difficult concepts and ideas. Now and again, a concept is stated without explanation, but references to the literature are provided. An elementary background in finite fields analysis and number theory at the level of [MC87], [LN97], [ST92], [WN03] etc, should be sufficient to read the material and will be assumed.

Due to its applications in cryptography, the amount of research in this subject has accelerated in the last few decades and there is a vast and growing literature on the subject. An effort was made to use only standard terms, notations, and symbols. This ensures a consistent and coherent interface among the various publications. The most current results are posted in online journals and preprint servers.

The earlier sections explain the standard notations and definitions of elliptic curves and the application of elliptic curves to cryptography. Some protocols and standards are discussed, but not extensively, see [NT99], [SECQ], [BL05] for detailed analysis. The hardware and software aspect of elliptic curves cryptography are treated in [HM04], [FT06], and the references therein. The remaining sections focus on the theory of the groups of points of elliptic curves, algebraic curves and related subjects.

# 2. The Discrete Logarithm Problems

There are a few open problems related to the discrete logarithms in cyclic groups. The problems emerged after the introduction of the Diffie-Hellman protocol in 1976. These open problems are stated below.

(i) *The Discrete Logarithm Problem*. Let $G = \{ g^r : 0 \leq r < n \}$ be a random cyclic group of $n$ elements, generated by $g > 1$. Given a primitive element $g$ and a random element $s = g^r \in G$, find the integer $r = \log(s)$.

(ii) *Computational DH Problem*. Given $n$, $g$, and a random pair of elements $g^r$ and $g^s$, find $g^{rs}$.

(iii) *Decision DH Problem*. Given $n$, $g$, and a random triple of elements $g^r$, $g^s$, $g^t$, determine whether or not $g^t = g^{rs}$.

(iv) *Indistinguishability Assumption*. A randomly generated triple $(g^r, g^s, g^{rs}) \in G \times G \times G$ is polynomial time indistinguishable from a randomly selected triple $(x, y, z) \in G \times G \times G$.

The calculation of the discrete logarithm on a random cyclic group and the related problems are considered intractable problems. The fastest algorithms run in exponential times $O(n^{1/2+\varepsilon})$, $\varepsilon > 0$. The running time complexity measures the number of group operations required to determine a single instance of the discrete logarithm. Over thirty years of extensive research on this problem (since the introduction of the DH protocol) has not changed the time complexity of the discrete logarithm in random cyclic groups.

Elliptic curves cryptography rests on the intractability of the discrete logarithm problem of the groups of points of elliptic curves. The fastest known discrete logarithm algorithms on the groups of points of random algebraic curves of small genus have exponential time complexities. Moreover, the groups of points of algebraic curves of small genus $g > 0$ such as elliptic curves offer significantly smaller public keys parameters at the same level of security as the established public keys systems, such as integers exponentiations. But, most importantly, there is no known subexponential algorithm to compute discrete logarithms in the groups of points of random algebraic curves of small genus.

The *genus* of a curve is a measure of the arithmetic and geometric complexity of the curve, [HN77]. The simplest curves, such as lines $ax + by + cz = 0$ and conics $ax^2 + bxy + cy^2 + dz^2 = 0$, have genus $g = 0$. These are followed by elliptic curves of genus $g = 1$, and hyperelliptic curves of small genus $g = 2, 3$, etc. The genus is essentially given by the expression $g = (d - 1)(d - 2)/2$, with $d$ being the maximal degree.

It should be observed that the discrete logarithms on some special classes of algebraic curves of small genus $g > 0$ do have algorithms of polynomial time or subexponential time complexities, see [MV93], [SV98], [SM99].

In contrast, there are subexponential algorithms to compute discrete logarithms in any finite field $\mathbf{F}_q$, see [ME97, p. 103], and there appears to be subexponential algorithms to compute discrete logarithms in the groups of points of every algebraic curves genus $g = 0$, and large $g \gg 4$. For example, the discrete logarithm in the groups of points of the conic $x^2 - y^2 = c$ over $\mathbf{F}_q$, $c$ constant, has the same complexity as the discrete logarithm in $\mathbf{F}_q$, refer to [LR00].

## 3. Definitions of Elliptic Curves

An elliptic curve over an arbitrary field $\mathbf{F}$ (or ring $\mathbf{R}$) can be described by a polynomial equation $f(x, y) \in \mathbf{F}[x,y]$ of maximal degree 3 or 4 in each variable. Various standard equations are employed in the theory of elliptic curves. The forms of these equations are contingent on the characteristic of the fields of definition. Elementary introductions to elliptic curves are given in [WN05], [MS93], [HM04], etc., and advanced coverage in [TE74], [CA66], [SN92], [FT06], etc.

*Weirstrauss Long Form*
The *Weirstrauss long form* of an elliptic curve in a field of any characteristic is defined by

$$E : y^2 + a_1xy + a_3y = x^3 + a_2x^2 + a_4x + a_6. \qquad (1)$$

The derivation of this equation (Riemann-Roch Theorem) is quite standard, see [CA66], [TE74], [HR87, p.66], [FT06], etc. This is an affine model with solutions in the affine plane $\mathbf{F} \times \mathbf{F}$. The set of $\mathbf{F}$-rational points is the set of solutions

$$E(\mathbf{F}) = \{ (x, y) : f(x, y) = 0 \text{ with } x, y \in \mathbf{F} \} \cup \{ O \}. \qquad (2)$$

The extra point represented by $O = (\infty, \infty)$ is called the point at infinity. The corresponding projective model is given by $z^3 f(x/z, y/z) = 0$. Specifically,

$$E : y^2z + a_1xyz + a_3yz^2 = x^3 + a_2x^2z + a_4xz^2 + a_6z^3. \qquad (3)$$

A solution of the projective model is a triple $(x, y, z) \in \mathbf{F} \times \mathbf{F} \times \mathbf{F}$, and the set of $\mathbf{F}$-rational points is the set of solutions $E(\mathbf{F}) = \{ (x, y, z) : f(x, y, z) = 0 \text{ with } x, y, z \in \mathbf{F} \} \cup \{ O \}$. Both the affine and projective equations are widely used in theory and practice.

*Tate Terms.* The Tate terms are derived from the coefficients of the Weirstrauss long form equation (1). These terms are used to compute various elliptic curves' parameters and to transform the long form equation to simpler forms.

$$b_2 = a_1^2 + 4a_2, \qquad b_4 = 2a_4 + a_1a_3, \qquad b_6 = a_3^2 + 4a_6, \qquad (4)$$
$$b_8 = a_1^2 a_6 + 4a_2 a_6 - a_1 a_3 a_4 + a_2 a_3^2 - a_4^2, \qquad 4b_8 = b_2 b_6 - b_4^2,$$
$$c_4 = b_2^2 - 24b_4, \qquad c_6 = -b_2^3 + 36 b_2 b_4 - 216 b_6.$$

The derivations of these relations are straight-forward simple changes of variables, see [TE74]. Another important invariant of a curve is derived from the roots of the defining polynomial $f(x) = (x - e_1)(x - e_2)(x - e_3)$ of the elliptic curve $E : y^2 = f(x)$.

***Definition 1.*** The *discriminant* of an elliptic curve is defined by the product of the squared differences $\Delta(E) = 16(e_1 - e_2)^2 (e_1 - e_3)^2 (e_2 - e_3)^2$ of the roots of the equation $f(x) = 0$. Expressed in terms of the coefficients, the discriminant has the form

$$\Delta(E) = -b_2^2 b_8 - 8b_4^3 - 27 b_6^2 + 9 b_2 b_4 b_6 = 12^{-3}(c_4^3 - c_6^2). \qquad (5)$$

An algebraic curve is called *nonsingular* if and only if it has distinct roots, or equivalently $\Delta(E) \ne 0$. Otherwise, the curve is called *singular*. The nonsingularity of an algebraic curve is also specified by the nonvanishing of the partial derivatives $\partial_x f(x, y) \ne 0$, and $\partial_y f(x, y) \ne 0$ at each **F**-rational points $(x, y) \in E(\mathbf{F})$, see [HR87].

***Definition 2.*** The *j-invariant* of an elliptic curve is defined by

$$j(E) = (b_2^2 - 24b_4)^3 / \Delta(E) = c_4^3 / \Delta(E) = 12^3 c_4^3 /(c_4^3 - c_6^2). \qquad (6)$$

***Example 3.*** The discriminant and *j*-invariant of the elliptic curve $E : y^2 + xy + 2y = x^3 + 3x^2 + 4x + 5$ are

$$\Delta(E) = -b_2^2 b_8 - 8b_4^3 - 27 b_6^2 + 9 b_2 b_4 b_6 = 2^2 \cdot 5^2 \cdot 139, \text{ and } j(E) = (b_2^2 - 24b_4)^3 / \Delta(E) = -2351^3 /(2^2 \cdot 5^2 \cdot 139),$$

where
$$b_2 = a_1^2 + 4a_2 = 7, \qquad b_4 = 2a_4 + a_1 a_3 = 10,$$
$$b_6 = a_3^2 + 4a_6 = 28, \qquad b_8 = a_1^2 a_6 + 4a_2 a_6 - a_1 a_3 a_4 + a_2 a_3^2 - a_4^2 = 44.$$

This example is intended to show how to use the Tate constants in the computation of the parameters. However, in most applications, almost all the coefficients $a_i$ are zero and the calculations are easier. A few techniques used to simplify the calculations and other aspects of elliptic curves are taken up next.

***Definition 4.*** Two elliptic curves, $E$ and $E'$, are *isomorphic* if there exists an *admissible transformation* and its inverse

$$(x, y) = (u^2 x' + r, u^3 y' + u^2 s x' + t) \text{ and } (x', y') = (u^{-2}(x - r), u^{-3}(y - sx + rs - t)), \qquad (7)$$

respectively that maps $E$ onto $E'$.

***Proposition 5.*** A pair of elliptic curves $E$ and $E'$ are called *isomorphic* over some finite extension of the field of definition **F** if and only if the *j*-invariants are equal: $j(E) = j(E')$.

In characteristic $p \ne 2, 3$, the admissible transformation

$$(x, y) \to \left( \frac{x - 3a_1^2 - 12a_2}{36}, \frac{y - 3a_1 x}{216} - \frac{a_1^3 + 4a_1 a_2 - 12a_3}{24} \right) \tag{8}$$

transforms the long form equation (1) into the short form $E': y^2 = x^3 - 27c_4 x - 54c_6$.

Due to this fact, it can be assumed that $a_1 = a_2 = a_3 = 0$, so that the long form reduces to the short form equation

$$E : y^2 = x^3 + a_4 x + a_6. \tag{9}$$

The discriminant and $j$-invariant of an elliptic curve defined by the short form equation are $\Delta(E) = -16(4a_4^3 + 27a_6^2)$, and $j(E) = 4 \cdot 12^3 a_4^3 / (4a_4^3 + 27a_6^2)$, respectively. The two special cases

| | | |
|---|---|---|
| $E : y^2 = x^3 + a_4 x$, | $\Delta(E) = -2^6 a_4^3$, | $j(E) = 1728$; |
| $E' : y^2 = x^3 + a_6$, | $\Delta(E') = -2^4 \cdot 3^3 a_6^2$, | $j(E') = 0$, |

play important roles in the theory and practice of elliptic curves.

***Example 6.*** This example introduces simple computation ideas used to determine the groups of points of elliptic curves, in particular for $E : y^2 = x^3 + 2$ over $\mathbf{F}_{13}$. Since its $j$-invariant and discriminant are $j(E) \equiv 0 \bmod 13$ $\Delta(E) = -16(27 \times 2^2) \equiv 1 \bmod 13$, respectively, this is a nonsingular curve over the finite field $\mathbf{F}_{13}$.

To determine the points on the curve, the values of the polynomial $f(x) = x^3 + 2$ are evaluated in a table. Those values that are square modulo 13 contribute one or two points to the group of points.

| X | $f(x) = x^3 + 2$ | X | $f(x) = x^3 + 2$ | x | $f(x) = x^3 + 2$ |
|---|---|---|---|---|---|
| −6 | 7 | −1 | 1 | 4 | 1 |
| −5 | 7 | 0 | 2 | 5 | 10 |
| −4 | 3 | 1 | 3 | 6 | 10 |
| −3 | 1 | 2 | 10 | | |
| −2 | 7 | 3 | 3 | | |

Note that the set of squares in $\mathbf{F}_{13}$ is

$$\begin{aligned} Q &= \{ x^2 : 0 \neq x \in \mathbf{F}_{13} \} \\ &= \{ (\pm 1)^2 = 1, (\pm 2)^2 = 4, (\pm 3)^2 = 9, (\pm 4)^2 = 3, (\pm 5)^2 = 12, (\pm 6)^2 = 10 \} \\ &= \{ 1, 3, 4, 9, 10, 12 \}. \end{aligned}$$

Now select the values that are squares in $\mathbf{F}_{13}$, for instance, $x = -4$ evaluates to $f(-4) = 3$, and $(\pm 4)^2 = 3$, so $x = -4$ contributes two points $(x, y) = (-4, \pm 4)$ to the set. The complete set of points is

$$E(\mathbf{F}_{13}) = \{ (\infty, \infty), (-4, \pm 4), (-3, \pm 1), (-1, \pm 1), (1, \pm 4), (2, \pm 6), (3, \pm 4), (4, \pm 1), (5, \pm 6), (6, \pm 6) \}.$$

The group has $19 = \#E(\mathbf{F}_{13})$ points including the identity point $(\infty, \infty)$, so the set $E(\mathbf{F}_{13})$ is isomorphic to a cyclic group, viz, $E(\mathbf{F}_{13}) \cong \mathbb{Z}_{19}$.

## 4. Elliptic Curves in Characteristic 2

There are two equations, (10) and (11), commonly used in elliptic curve cryptography in characteristic 2. The determination of these equations from the long form equation (1) is considered now.

**Case 1.** Let $a_1 = 0$, $a_3 \neq 0$. Then the Tate constants are:

$b_2 = a_1^2 + 4a_2 = 0$, $\qquad b_4 = 2a_4 + a_1 a_3 = 0$,

$b_6 = a_3^2 + 4a_6 = a_3^2$, $\qquad b_8 = a_1^2 a_6 + 4a_2 a_6 - a_1 a_3 a_4 + a_2 a_3^2 - a_4^2 = a_2 a_3^2$.

The discriminant reduces to $\Delta(E) = -b_2^2 b_8 - 8b_4^3 - 27b_6^2 + 9b_2 b_4 b_6 = a_3^4$, and the $j$-invariant reduces to $j(E) = (b_2^2 - 24b_4)^3 / \Delta(E) = 0$.

Now for any $(a_3, a_4) \neq (0, 0)$, or $(a_3, a_6) \neq (0, 0)$, the discriminant $\Delta(E) \neq 0$ does not vanish, and $j(E) = 0$. Thus, the equation $E : y^2 + a_3 y = x^3 + a_2 x^2 + a_4 x + a_6$ defines a nonsingular algebraic curve.

The (admissible) linear change of variables $(x, y) \to (x + a_2, y)$ maps the previous equation to the more common form

$$E' : y^2 + cy = x^3 + ax + b \qquad (10)$$

encountered in the literature ($c = a_3$, $a = a_4$, $b = a_6$). The corresponding discriminant and $j$-invariant are $\Delta(E') = c^4$, and $j(E') = 0$ respectively. Every supersingular elliptic curve in characteristic 2 has an equation of this form. Moreover, there are no nontrivial points of even order. Thus, supersingular elliptic curves over $\mathbf{F}_{2^n}$ have trivial torsion $E[2^k] = \{ O \}$.

There are three isomorphism classes of supersingular curves over $\mathbf{F}_{2^n}$, if $n$ is odd, and seven isomorphism classes if $n$ is even.

**Case 2.** Let $a_1 \neq 0$, $a_3 = 0$. Then the Tate constants are:

$b_2 = a_1^2 + 4a_2 = a_1^2$, $\qquad b_4 = 2a_4 + a_1 a_3 = 0$,

$b_6 = a_3^2 + 4a_6 = 0$, $\qquad b_8 = a_1^2 a_6 + 4a_2 a_6 - a_1 a_3 a_4 + a_2 a_3^2 - a_4^2 = a_1^2 a_6$.

The discriminant reduces to $\Delta(E) = -b_2^2 b_8 - 8b_4^3 - 27b_6^2 + 9b_2 b_4 b_6 = a_1^2(a_1^2 a_6 + a_4^2)$, and the $j$-invariant reduces to $j(E) = (b_2^2 - 24b_4)^3 / \Delta(E) = a_1^{10}(a_1^2 a_6 + a_4^2)^{-1}$.

Now for any $(a_1, a_4) \neq (0, 0)$, or $(a_1, a_6) \neq (0, 0)$, the discriminant $\Delta(E) \neq 0$ does not vanish, and $j(E) \neq 0$. Thus, the equation $E : y^2 + a_1 xy = x^3 + a_2 x^2 + a_4 x + a_6$ defines a nonsingular algebraic curve.

The (admissible) linear change of variables $(x, y) \to (a_1^2 x + a_3 a_1^{-1}, a_1^3 y + (a_1^2 a_4 + a_4) a_1^{-3})$ maps the previous equation to the more common form

$$E'' : y^2 + cxy = x^3 + ax^2 + b \qquad (11)$$

encountered in the literature ($c = a_1$, $a = a_2$, $b = a_6$). The corresponding discriminant and $j$-invariant are $\Delta(E'') = bc^6$, and $j(E'') = c^8/b$ respectively. Every nonsupersingular elliptic curve in characteristic 2 has an equation of this form. Moreover, there is a point $P = (x, y)$ of order 2 if $a_2 \neq 0$, otherwise $a_2 = 0$ and there is a point of order 4. Thus nonsupersingular elliptic curves over $\mathbf{F}_{2^n}$ have nontrivial torsion $E[2^k] = \mathbb{Z}_{2^k}$, $k \geq 1$.

***Example 7.*** (Hand calculations) As the coefficients $a_1$, $a_3$, $a_2$, $a_4$, $a_6$ of (1) vary over $\mathbf{F}_2$, 32 distinct equations are generated. Twenty of these equations are nonsingular, namely, $\Delta(E) \neq 0$. Twelve of the nonsingular equations have $j$-invariant $j(E) = 1$, and eight have $j(E) = 0$. One curve from each of the five isomorphic classes is listed below.

The 12 nonsingular (ordinary) curves of $j$-invariant $j(E) = 1$ consist of two isomorphic classes.

$E_1 : y^2 + xy = x^3 + 1$, $\qquad\qquad E_1(\mathbf{F}_2) = \{ O, (0, 1), (1, 0), (1, 1) \} \cong \mathbb{Z}_4$,

$E_2 : y^2 + xy = x^3 + x^2 + 1$, $\qquad\qquad E_2(\mathbf{F}_2) = \{ O, (0, 1) \} \cong \mathbb{Z}_2$,

The specific class a curve belongs to is determined by the change of variables $(x, y) \to (x, x + y)$.
The 8 nonsingular curves of $j$-invariant $j(E) = 0$ consist of three isomorphic classes.

$E_3 : y^2 + y = x^3 + x$, $\qquad\qquad E_3(\mathbf{F}_2) = \{ O, (0, 0), (0, 1), (1, 0), (1, 1) \} \cong \mathbb{Z}_5$,

$E_4 : y^2 + y = x^3 + x + 1$ $\qquad\qquad E_4(\mathbf{F}_2) = \{ O \} \cong \mathbb{Z}_1$,

$E_5 : y^2 + y = x^3$, $\qquad\qquad E_5(\mathbf{F}_2) = \{ O, (0, 0), (0, 1) \} \cong \mathbb{Z}_3$,

The specific class a curve belong to is determined by the change of variables $(x, y) \to (x + 1, x + y)$.

Other routine information about these five curves are tabulated below, these are computed directly from the previous calculations. For example, the characteristic polynomial $f_\pi(x) = x^2 - ax + q$ of $E_1$ is obtained from the values $q = 2$, and $a = q + 1 - \#E_1(\mathbf{F}_2) = q + 1 - 4$. This elliptic curve has complex multiplication by $\pi = (-1 + \sqrt{-7})/2$ because this complex number is a root $f_\pi(x)$. The number $\pi$ is an eigenvalue of the Frobenious map $(x, y) \to (x^q, y^q)$. Complex multiplication finds applications in point counting methods and scalar multiplications, see [SA00].

| Curve | $\pi$ | $f_\pi(x)$ | Group |
|---|---|---|---|
| $E_1$ | $\pi = (-1 + \sqrt{-7})/2$ | $x^2 + x + 2$ | $E_1(\mathbf{F}_2) \cong \mathbb{Z}_4$ |
| $E_2$ | $\pi = (1 + \sqrt{-7})/2$ | $x^2 - x + 2$ | $E_2(\mathbf{F}_2) \cong \mathbb{Z}_2$ |
| $E_3$ | $\pi = -1 + \sqrt{-1}$ | $x^2 + 2x + 2$ | $E_3(\mathbf{F}_2) \cong \mathbb{Z}_5$ |
| $E_4$ | $\pi = 1 + \sqrt{-1}$ | $x^2 - 2x + 2$ | $E_4(\mathbf{F}_2) \cong \mathbb{Z}_1$ |
| $E_5$ | $\pi = \sqrt{-2}$ | $x^2 + 2$ | $E_5(\mathbf{F}_2) \cong \mathbb{Z}_3$ |

***Proposition 8.*** There are five inequivalent isomorphism classes of elliptic curves over $\mathbf{F}_2$: 2 classes $E_1$ and $E_2$ of $j$-invariants $j(E_i) = 1$, and three inequivalent isomorphism classes $E_3$, $E_4$ and $E_5$ of $j$-invariants $j(E_i) = 0$.

Proof: Brute force calculations as above in Example 7. ∎

## 5. Elliptic Curves in Characteristic 3
In characteristic $p > 2$, the linear transformation $(x, y) = (x, y - (a_1 x + a_3)/2)$ transforms the long form equation (1) into the short form

$$E' : y^2 = x^3 + (b_2 x^2 + 2b_4 x + b_6)/4. \qquad (12)$$

Since $4 \equiv 1 \bmod 3$, the Tate constants are:

$$b_2 = a_1^2 + 4a_2 = a_2, \qquad\qquad b_4 = 2a_4 + a_1a_3 = 2a_4,$$
$$b_6 = a_3^2 + 4a_6 = a_6, \qquad\qquad b_8 = a_1^2 a_6 + 4a_2 a_6 - a_1 a_3 a_4 + a_2 a_3^2 - a_4^2 = 2a_2 a_6 - a_4^2.$$

The discriminant reduces to $\Delta(E) = a_2^2 a_4^2 - a_2^3 a_6 - a_4^3$, and the j-invariant reduces to $j(E) = a_2^6 / (a_2^2 a_4^2 - a_2^3 a_6 - a_4^3)$. There are two important cases:

**Case 1.** $a_1 = a_2 = a_3 = 0$, and $(a_4, a_6) \neq (0, 0)$. For these parameters the discriminant $\Delta(E) \neq 0$ does not vanish, and $j(E) = 0$. Thus, the equation $E : y^2 = x^3 + a_4 x + a_6$ defines a nonsingular algebraic curve.

**Case 2.** $a_1 = a_3 = 0$, $a_2 \neq 0$, and $(a_4, a_6) \neq (0, 0)$. For these parameters the discriminant $\Delta(E) \neq 0$ does not vanish, and $j(E) \neq 0$. Thus, the equation $E : y^2 = x^3 + a_2 x^2 + a_4 x + a_6$ defines a nonsingular algebraic curve.
For example, the curve $E : y^2 = x^3 + x^2 + a_6$ has the j-invariant $j(E) = -1/a_6$, so every value of $0 \neq j(E) \in \mathbf{F}_3$ is possible.

The following elliptic curves are representatives:

| | | |
|---|---|---|
| $E_1 : y^2 = x^3 + x,$ | $E_2 : y^2 = x^3 - x,$ | j-invariant $j(E) = 0,$ |
| $E_3 : y^2 = x^3 - x + 1,$ | $E_4 : y^2 = x^3 - x - 1,$ | j-invariant $j(E) = 0,$ |
| $E_5 : y^2 = x^3 + x^2 - 1$ | $E_6 : y^2 = x^3 - x^2 + 1,$ | j-invariant $j(E) = 1$ |
| $E_7 : y^2 = x^3 + x^2 + 1$ | $E_8 : y^2 = x^3 - x^2 - 1,$ | j-invariant $j(E) = -1.$ |

***Proposition 9.*** There are eight isomorphic inequivalent classes of elliptic curves over $\mathbf{F}_3$: There are four inequivalent classes of j-invariants $j(E_i) = 0$, there are two inequivalent classes of j-invariants $j(E_i) = 1$, and there are two inequivalent classes of j-invariants $j(E_i) = -1$.

Proof: Try computing the groups of points $E(\mathbf{F}_3)$ of the 27 elliptic curves $E : y^2 = x^3 + ax^2 + bx + c$, with $a, b, c \in \mathbf{F}_3$ and linear maps between them, everything by hand. ∎

## 6. Standard Curves of Prescribed j-Invariants
The j-invariant determines a unique curve up to an isomorphism. This parameter classifies the representatives of elliptic curves of a given value of j. Standard curves with prescribed j-invariants in characteristic char $p > 0$ are tabulated below.

| Char $p$ | Discriminant | j-invariant | Curve |
|---|---|---|---|
| $p = 2$ | $\Delta = a_3^4$ | $j = 0$ | $y^2 + a_3 y = x^3 + a_4 x + a_6$ |
| $p = 2$ | $\Delta = a_6$ | $j = a_6^{-1}$ | $y^2 + a_1 xy + = x^3 + a_2 x^2 + a_6$ |
| $p = 3$ | $\Delta = -a_4^3$ | $j = 0$ | $y^2 = x^3 + a_4 x + a_6$ |
| $p = 3$ | $\Delta = -a_2^3 a_6$ | $j = -a_2^3 a_6^{-1}$ | $y^2 = x^3 + a_2 x^2 + a_6$ |
| $p \neq 2, 3$ | $\Delta = -2^4 \cdot 3^3 a_6^2$ | $j = 0$ | $y^2 = x^3 + a_6$ |
| $p \neq 2, 3$ | $\Delta = -2^6 a_4^3$ | $j = 1728$ | $y^2 = x^3 + a_4 x$ |
| $p \neq 2, 3$ | $\Delta = -16[4(-a)^3 + 27(2a)^2]$ | $(E) = j \neq 0, 1728$ | $E : y^2 = x^3 - ax \pm 2a,$ $a = 27j/(j - 12^3)$ |

***Proposition 10.*** There are precisely $q^2 - q$ nonsingular elliptic curves $E : y^2 = x^3 + ax + b$ over $\mathbf{F}_q$, $q \neq 2^v, 3^v$, $v \geq 1$.

Proof: There are $q^2$ distinct curves in characteristic $p > 3$, one for each pair $(a, b)$. The curve is nonsingular if and only if the discriminant $-16(4a^3 + 27b^2) \neq 0$. On the other hand, the equation $4a^3 + 27b^2 = 0$ has precisely $q$ solutions, namely, $a = -3u^2, b = 2u^3, u \in \mathbf{F}_q$. ∎

***Proposition 11.*** Let $E : y^2 = x^3 + ax + b$ and $E' : y^2 = x^3 + \bar{a}x + \bar{b}$. Then
(i) The curves $E$ is isomorphic to $E'$ if and only if there exists $u \neq 0$ such that $a = u^4\bar{a}, b = u^6\bar{b}$.
(ii) There are $2q + 6, 2q + 2, 2q + 4, 2q$ isomorphic classes of elliptic curves over $\mathbf{F}_q$ for $q \equiv 1, 5, 7, 11$ mod 12, respectively.
(iii) If $q = 2^n$, and $n$ is odd, then there are $2^{n+1} + 1$ isomorphic classes of elliptic curves over $\mathbf{F}_q$, otherwise, $n$ is even and there are $2^{n+1} + 5$ isomorphic classes.
(iv) If $q = 3^n$, and $n$ is odd, then there are $2(3^n + 1)$ isomorphic classes of elliptic curves over $\mathbf{F}_q$, otherwise, $n$ is even and there are $2(3^n + 1)$ isomorphic classes.

The proof of the first statement follows from the definition, and the proof of the other statements appear in [SF87]. The isomorphism classes of curves are represented by $(q - 1)/w$ representatives of distinct curves $E : y^2 = x^3 + ax + b$ as $(a, b)$ ranges over the set of ordered pairs such that $4a^3 + 27b^2 \neq 0$. Here $w$ is the number of automorphisms of the curve, that is, $w = 2$ if the $j$-invariant $j \neq 0, 1728$, and $w = 4, 6$ if the $j$-invariant $j = 0, 1728$. To see this, note that for $a^3/b^2 = -27/4$ the curve is singular, but for $a^3/b^2 = c \neq 0, -27/4$, there are two isomorphic classes of curves such that $a = u^4a', b = u^6b'$ with $u \neq 0$. For the other cases $a = 0, b \neq 0$, and $a \neq 0, b = 0$, there are four and six isomorphic classes of curves, respectively.

It should be observed that an isomorphism between two curves implies that $E(\mathbf{F}_q) \cong E'(\mathbf{F}_q)$, but not conversely. In fact, $E(\mathbf{F}_q) \cong E'(\mathbf{F}_q)$ does not imply that the elliptic curves are isomorphic. For example, take $E : y^2 = x^3 + 1$ and $E' : y^2 = x^3 + b$ over $\mathbf{F}_q$ such that $b \neq c^6$ and $q \equiv 2$ mod 3. Then both $E(\mathbf{F}_q)$ and $E'(\mathbf{F}_q)$ have the same order $q + 1$, but the elliptic curves are not isomorphic since there is no admissible map between the curves.

***Example 12.*** There are $q^2 - q = 110$ nonsingular elliptic curves with coefficients in $\mathbf{F}_{11}$. These are classified into $2q = 22$ nonisomorphic classes of $(q - 1)/2 = 5$ curves each. By mean of the map $u \to (a, b) = (u^4, u^6)$, where $u = \pm 1, \pm 2, \pm 3, \pm 4, \pm 5$. In particular, for $(a, b) = (1, 1)$, the following curves are in the same isomorphic class:

$E_1 : y^2 = x^3 + x + 1,$       $E_2 : y^2 = x^3 + 5x - 2,$       $E_3 : y^2 = x^3 + 4x + 3,$
$E_4 : y^2 = x^3 + 3x + 4,$      $E_5 : y^2 = x^3 - 2x + 5.$

Note that $\Delta(E_1) \equiv -1$ mod 11 $\not\equiv \Delta(E_2) \not\equiv \Delta(E_3) \not\equiv \Delta(E_4) \not\equiv \Delta(E_5)$, but $j(E_1) \equiv -2$ mod 11 $\equiv j(E_2) \equiv j(E_3) \equiv j(E_4) \equiv j(E_5)$.

An *isogeny* $\phi : E \to E'$ between two elliptic curves is a rational map given by

$$\phi(x, y) = \left( \frac{r(x)}{s(x)}, cy \frac{d}{dx} \frac{r(x)}{s(x)} \right), \tag{13}$$

where $r(x), s(x) \in \mathbf{F}_q[x]$ are polynomials, and $c^{12} = \Delta(E')/\Delta(E)$, see [SN2, p.183] for more details. Isogenies come in matched pairs $\phi, \hat{\phi}$ such that $\phi \circ \hat{\phi}(P) = [m]P$, where $m \geq 1$ is the degree of the map.

***Proposition 13.*** Let $E$ and $E'$ be a pair of elliptic curves over the finite field $\mathbf{F}_q$. Then the curves are *isogenous* if and only if $\#E(\mathbf{F}_q) = \#E'(\mathbf{F}_q)$.

The groups $E(\mathbf{F}_q)$ and $E'(\mathbf{F}_q)$ of a pair of isogenous curves can have the same cardinality but can have different structures, i. e., one cyclic and the other noncyclic, or different $j$-invariants.

## 7. Endomorphisms

Given a set $S$, an *endomorphism* on $S$ is a map from a set to itself. The set of all endomorphisms on the set of points $E(\mathbf{F}_q) = \{ \phi : E \to E \}$ of an elliptic curve to itself is denoted by $\mathrm{End}(E)$. An endomorphism is a rational map $\phi(P) = (f(P), g(P))$, where $f(x,y)$, $g(x,y)$ are rational functions, and $P \in E(\mathbf{F}_q)$. The set $\mathrm{End}(E)$ form a ring:

(i) $\phi(O) = O$,     identity,
(ii) $\phi(P + Q) = \phi(P) + \phi(Q)$,     additive,
(iii) $\phi(nP) = n\phi(P)$,     isobaric.

The set of endomorphisms on the algebraic curve have the following properties:
(1) It is a well-defined map at every point of the algebraic curve $C$.
(2) It is a homomorphism with respect to the group law on the algebraic curve $C$.

To form a ring of maps, the set of endomorphisms must satisfy the second property above. Otherwise the distributive law on the set of points of the algebraic curve $C$ would fail.

An endomorphism $\pi \in \mathrm{End}(E)$ of a nonCM elliptic curve has a linear or quadratic characteristic polynomial $f_\pi(x) = x^2 - Tr(\pi)x + N(\pi)$, where $Tr(\pi) = \pi + \bar{\pi}$, and $N(\pi) = \pi\bar{\pi}$ are the trace and norm.

Two important endomorphisms of elliptic curves are the followings.
(i) The scalar point multiplication is the rational map $P \to [n]P$, which is an endomorphism of degree $n^2$.

(ii) The Frobenius map $\sigma : \overline{\mathbf{F}}_q \to \overline{\mathbf{F}}_q$ defined by $\sigma(x) = x^q$. This map fixes the elements of the finite field $\mathbf{F}_q$. Here, $\overline{\mathbf{F}}_q = \bigcup_{n=1}^{\infty} \mathbf{F}_{p^n}$ is the *algebraic closure* of the finite field $\mathbf{F}_q$. The algebraic closure contains all the roots of every algebraic equation over $\mathbf{F}_q$. The Frobenius map has a natural extension to the product space $\sigma : \overline{\mathbf{F}}_q \times \overline{\mathbf{F}}_q \times \cdots \times \overline{\mathbf{F}}_q \to \overline{\mathbf{F}}_q \times \overline{\mathbf{F}}_q \times \cdots \times \overline{\mathbf{F}}_q$ by $(x_0, x_1, ..., x_n) \to (x_0^q, x_1^q, ..., x_n^q)$. The degree of the this map is $\deg(\sigma) = q$, and the trace of this map is given by the sum

$$Tr(\sigma) = \sum_{i=1}^{2g} \alpha_i = \sum_{i=1}^{g} \alpha_i + \bar{\alpha}_i \tag{14}$$

of the roots $\alpha_1, \alpha_2, ..., \alpha_g, \alpha_{g+1} = \bar{\alpha}_1, \alpha_{g+2} = \bar{\alpha}_2, ..., \alpha_{2g} = \bar{\alpha}_g$ of the characteristic polynomial of $\sigma$.

The map $\sigma$ corresponds to complex multiplication by $\sigma = (t + u\sqrt{d})/2 \in \mathrm{End}(E)$, where $f_\sigma(x) = x^2 - tx + q$ is the characteristic polynomial, $t = Tr(\sigma) = \sigma + \bar{\sigma}$ is the trace, and $q = N(\sigma) = \sigma\bar{\sigma}$ is the norm of the map.

*Proposition* 14. The ring of endomorphisms of an elliptic curve is one of the following:
(i) $\mathrm{End}(E) \cong \mathbb{Z}$.     (ii) $\mathrm{End}(E) \cong$ Nonmaximal Order in $\mathbb{Q}(\sqrt{d})$.
(iii) $\mathrm{End}(E) \cong$ Maximal Order in a Quaternion Ring $\mathbb{Q}(i,j,k)$.     (iv) $\mathrm{End}(E) \cong$ Maximal Order in $\mathbb{Q}(\sqrt{d})$.

***Definition* 15.**  An elliptic curve is said to have *complex multiplication* (denoted by CM) if $\mathbb{Z} \subset \text{End}(E)$, but $\text{End}(E) \not\subset \mathbb{Z}$.

Since a map $\phi \in \text{End}(E)$ can be identified with an element $\pi \in \mathbb{Z}$, or $\mathbb{Q}(\sqrt{d}\,)$, or $\mathbb{Q}(i,j,k)$, the action of $\phi$ on $E$ is actually scalar multiplication $\phi(P) = [\pi]P$.

For example, in finite fields every elliptic curve has complex multiplication since the Frobenius map $\sigma$ of a nonsingular elliptic curve has the representation $\sigma(P) = [(t + \sqrt{d}\,)/2]P$, with $t = q + 1 - N$, and $d < 0$. Specifically, for ordinary elliptic curves, the Frobenius map $\sigma(P) = (x^q, x^q) \neq [n]P$ is not multiplication by an integer, so the group $\text{End}(E)$ contains the set of integers. On the other hand, Frobenius map $\sigma$ of a supersingular elliptic curve has the representation $\sigma(P) = (x^{p^2}, y^{p^2}) = [\pm p]P$ is multiplication by an integer, nevertheless, the group $\text{End}(E)$ is the maximal order in a quaternion ring.

The best known elliptic curves over the integers $\mathbb{Z}$ with CM are the two special cases:

(i) $E: y^2 = x^3 + a_4 x$, $\quad\quad\quad \Delta(E) = -2^6 a_4^3$, $\quad\quad\quad j(E) = 1728;$ $\quad\quad\quad \text{End}(E) \cong \mathbb{Z}[\sqrt{-1}\,].$

In this case, there are four elliptic curves with the same *j*-invariant and endomorphism group. Here the complex multiplication is defined as

$$\sigma(P) = [(t + \sqrt{-1}\,)/2]P = (-x, y\sqrt{-1}\,), \text{ where } P = (x, y).$$

(ii) $E: y^2 = x^3 + a_6$, $\quad\quad\quad \Delta(E) = -2^4 \cdot 3^3 a_6^2$, $\quad\quad\quad j(E) = 0;$ $\quad\quad\quad \text{End}(E) \cong \mathbb{Z}[\sqrt{-3}\,].$

In this case, there are six elliptic curves over the integers $\mathbb{Z}$ with the same *j*-invariant and endomorphism group. Here, the complex multiplication is defined as

$$\sigma(P) = [(-1 + \sqrt{-3}\,)/2]P = (x(-1 + \sqrt{-3}\,)/2, y).$$

The theory of complex multiplication has prominent applications in elliptic curve cryptography.

## 8. Automorphisms

An endomorphism that leaves the elliptic curve equation unchanged is called an *automorphism*. The group of automorphisms $\text{Aut } E = \{\sigma : E \to E\}$ of an elliptic curve has cardinalities of 2, 4, 6, 12 or 24 depending on the characteristic of the field of definition and the *j*-invariant. For most curves the group has 2 maps. The first three groups consist of the following maps:

$$(x, y) \to (x, -y) \quad\quad\quad (x, y) \to (\pm x, \pm iy) \quad\quad\quad (x, y) \to (\omega x, \pm y), \quad\quad\quad (15)$$

where $i^2 = -1$, and $\omega^3 = 1$, respectively.

***Proposition* 16.**  Given an elliptic curve over a finite, the following statements hold:
(i) If $q \neq 2^v, 3^v$, and $j(E) \neq 0, 1728$, $(a \neq 0, b \neq 0)$, then $\text{Aut } E \cong \{1, -1\}$.
(ii) If $q \neq 2^v, 3^v$, and $j(E) = 1728$, $(a \neq 0, b = 0)$, and $\mathbf{F}_q$ contains an element $\omega$ of order 4, then $\text{Aut } E \cong \mu_4$.
(iii) If $q \neq 2^v, 3^v$, and $j(E) = 0$, $(a = 0, b \neq 0)$, and $\mathbf{F}_q$ contains an element $\omega$ of order 6, then group of automorphisms $\text{Aut } E \cong \mu_6$.
(iv) If $q = 3^v$, and $j(E) = 0$, and $\mathbf{F}_q$ contains elements $\omega$, $\xi$ of orders 4 and 3, then group of automorphisms $\text{Aut } E$ is a noncommutative group that contains $\mu_4 \times \mu_3$.

(*v*) If $q = 2^v$, and $j(E) = 0$, and $\mathbf{F}_q$ contains elements $\omega$, $\xi$ of orders 6 and 4, then group of automorphisms Aut $E$ is a noncommutative group that contains $\mu_6 \times \mu_4$.

The notation $\mu_n = \{ \omega \in \mathbf{F} : \omega^n = 1 \}$ denotes the group of *n*th roots of unity in an arbitrary ring $\mathbf{R}$ or field $\mathbf{F}$, and $\mu_m \times \mu_n$ denotes the cross product of two such groups. For example, $\mu_4 = \{ 1, \omega, \omega^2, \omega^3 \}$, and $\mu_3 \times \mu_4 \cong \{ 1, \omega, \omega^2 \} \times \{ 1, \xi, \xi^2, \xi^3 \}$ in field (or ring) that contains these elements.

Since most ordered pairs are of the shape $(a \neq 0, b \neq 0)$, a random elliptic curve in characteristic $q \neq 2^v, 3^v$ has a very high probability of having just two automorphisms, see case (i). In contrast, the ordered pairs $(a \neq 0, b = 0)$, and $(a = 0, b \neq 0)$ are very rare, so a random elliptic curve has a very low probability of having four or six automorphisms, see cases (ii) and (iii).

## 9. Types of Singularities of Elliptic Curves in Characteristic $p > 3$

The discriminants of elliptic curves classify the curves into two types: nonsingular and singular. Nonsingular curves come in one type, and the singular curves are classified into three types, see [IC, p. 135] for more details.

1. Nonsingular curves $E : y^2 = x^3 + ax + b$, and discriminant $\Delta = -16(4a^3 + 27b^2) \neq 0$, trace $|a_p| < 2p^{1/2}$.
2. Singular curves $E : y^2 = x^3 + x^2$ with a node, rational slope, and discriminant $\Delta = -16(4a^3 + 27b^2) = 0$, $a \neq 0$, and trace $a_p = 1$.

3. Singular curves $E : y^2 = x^3 + x^2$ with a node, irrational slope, discriminant $\Delta = -16(4a^3 + 27b^2) = 0$, $a \neq 0$, and trace $a_p = -1$.

4. Singular curves $E : y^2 = x^3$ with a cusp, and discriminant $\Delta = -16(4a^3 + 27b^2) = 0$, $a = 0$, and trace $a_p = 0$.

The relationship between singularity and the cardinality of the group of points $E(\mathbf{F}_q)$ of an elliptic curve is encapsulated in the characteristic polynomial as follows: $\#E(\mathbf{F}_q) = qf_E(q^{-1})$, where

$$f_E(x) = \begin{cases} (1 - q^{1/2}x)(1 + q^{1/2}x) & \text{nonsingular over } \mathbf{F}_q, \\ 1 - x & \text{if the slope is rational over } \mathbf{F}_q, \\ 1 + x & \text{if the slope is irrational over } \mathbf{F}_q, \\ 1 & \text{if there is a cusp,} \end{cases} \tag{16}$$

see [TE74].

## 10. Supersingular Elliptic Curves

Supersingular elliptic curves are a special class of curves that have exceptional groups of points and other special structures. The property of supersingularity turns up in many different forms and there are many different ways of describing it. In characteristic $p = 2$ and 3 supersingular elliptic curves have *j*-invariants $j = 0$, but in characteristic $p \neq 2$ and 3 supersingular elliptic curves can have any *j*-invariant values.

**Definition 17.** Let $q = p^v$, $v > 0$. An elliptic curve is called *supersingular* if the following hold.
(i) $\#E(\mathbf{F}_q) = q + 1 - t \equiv q + 1 \bmod p$.
(ii) The $\pi^n$ is an integer for all even $n$ or all odd $n$, but not both, and $q = \pi\bar{\pi}$.
(iii) The $j$-invariant is in a quadratic extension of $\mathbf{F}_q$.
(iv) The endomorphism ring is a quaternion $\mathbb{Q}(i, j, k)$, with $i^2 = -1$, $j^2 = p$, $ij = -ji$,

An elliptic curve in characteristic $p > 0$ is supersingular if the torsion group $E[p^n] = \{\,O\,\}$ for all $n \geq 1$. Otherwise $E[p^n] \neq \{\,O\,\}$ and the curve is not supersingular.

**Example 18.** The two special cases of elliptic curves are supersingular at the following primes:
(i) $E : y^2 = x^3 + ax$, with $a \neq 0$, has $\#E(\mathbf{F}_q) = p^n + 1$ points for all $p \equiv 3 \bmod 4$, and all odd $n \geq 1$.
(ii) $E : y^2 = x^3 + b$, with $b \neq 0$, has $\#E(\mathbf{F}_q) = p^n + 1$ points for all $p \equiv 2 \bmod 3$, and all odd $n \geq 1$. If it has no roots over $\mathbf{F}_q$, the torsion $E[p] = \{\,O\,\} = E(\mathbb{Q})_{\text{tor}}$.

## 11. Legendre Form
An equation in *Legendre form* in characteristic $p > 2$ is described by

$$E_\alpha : y^2 = x(x - 1)(x - \alpha), \quad \alpha \neq 0, 1. \tag{17}$$

This curve satisfies the following:
(i) There are several points of order two.
$P = (0, 0)$,   $P = (\alpha, 0)$,   $P = (1, 0)$.
(ii) The discriminant and $j$-invariant are given by $\Delta(E_\alpha) = -16\alpha^2(\alpha - 1)^2$, and $j(E_\alpha) = 2^8(\alpha^2 - \alpha + 1)^3(\alpha - 1)^{-2}\alpha^{-1}$ respectively.
(iii) Two curves in $E_\alpha$ and $E_\beta$ are isomorphic if and only if

$$\beta \in \{\alpha, \alpha^{-1}, 1 - \alpha, (1 - \alpha)^{-1}, \alpha(\alpha - 1)^{-1}, (\alpha - 1)\alpha^{-1}\}. \tag{18}$$

The legendre form is a collection of curves with a torsion subgroup of order 4.

**Proposition 19.** An elliptic curve over $\mathbf{F}_q$ in characteristic $p > 2$ is isomorphic to a Legendre form over some extension of $\mathbf{F}_q$.

The proof is a straight forward application of a few linear changes of variables, see [SN92, p. 54].

## 12. Edwards Form
The Edwards form is a representation of an elliptic curve of the form

$$E_{c,d} : x^2 + y^2 = c^2(1 + dx^2 y^2), \tag{19}$$

where $c, d \neq 0$ and $cd^4 \neq 1$. There are several distinguished points.
(i) $P = (0, c) = \infty$, the point at infinity.     (ii) $-P = (-x, y)$, the additive inverse.
(ii) $P = (0, -c)$ has order 2.     (iii) $P = (-c, 0)$ and $(c, 0)$ have order 4.

The advantages of the Edwards form include speed and a simple uniform addition formula, see [BT07]. For a quadratic nonresidue $d \neq 0$, and a pair of points $P = (x_1, y_1)$ and $Q = (x_1, y_1)$, the addition formula is

$$P + Q = \left( \frac{x_1 y_2 + x_2 y_1}{c(1 + d x_1 x_2 y_1 y_2)}, \frac{y_1 y_2 - x_1 x_2}{c(1 - d x_1 x_2 y_1 y_2)} \right). \tag{20}$$

Unlike the standard addition formula, see Section 14, this formula holds for all points such that $x_1 x_2 y_1 y_2 \neq \pm 1$, for instance, it holds for $P = Q$.

The change of variables $w = y(1 - cx^2)$ maps the curve $E_{c,1} : x^2 + y^2 = c^2(1 + x^2 y^2)$ to the curve $E_c : y^2 = (c^2 - x^2)(1 - c^2 x^2)$.

**Lemma 20.** Let $b, c \neq 0$. The curves $E_{c,1} : x^2 + y^2 = c^2(1 + x^2 y^2)$ and $E_{b,1} : x^2 + y^2 = b^2(1 + x^2 y^2)$ are equivalent over $\mathbf{F}_q$ or $\mathbf{F}_q(i)$, if $b \neq 0$ is one of the values

$$ci^k, \quad c^{-1}i^k, \quad (c-1)(c+1)^{-1}i^k, \quad (c+1)(c-1)^{-1}i^k, \quad (c-i)(c+i)^{-1}i^k, \quad (c+i)(c-i)^{-1}i^k,$$

where $k = 0, 1, 2, 3$.

## 13. The Quadratic Twist Of A Curve

**Definition 21.** (i) Let $E : y^2 = x^3 + ax^2 + bx + c$ be a curve over $\mathbf{F}_q$, $q \neq 2^e$, and let $u > 0$ be a quadratic nonresidue. The quadratic twist of the curve $E$ is the associated curve $E' : y^2 = x^3 + aux^2 + bu^2 x + cu^3$.
(ii) Let $E : y^2 + (a_1 x + a_3) y = x^3 + a_2 x^2 + a_4 x + a_6$ be a curve over $\mathbf{F}_q$, $q = 2^e$, and let $\gamma$ be an element of trace $\mathrm{Tr}(\gamma) = 1$ in $\mathbf{F}_q$. The $\gamma$-twist of the curve $E$ is the associated curve $E' : y^2 = (a_1 x + a_3) y = x^3 + (a_2 + \gamma a_1^2) x^2 + a_4 x + a_6 + \gamma a_3^2$.

In characteristic $p > 2$, the $(q - 1)/2$ distinct twists of a curve are equivalent, and every twisted curve $E'$ is isomorphic to the curve $E'' : uy^2 = x^3 + ax^2 + bx + c$, $0 \neq u$ a nonsquare. This is obtained via the admissible change of variables $x \to ux$, $y \to u^2 y$, which maps $E'$ to $E''$. In contrast, in characteristic $p = 2$, there are $q/2$ distinct twists of a curve.

**Example 22.** If $\gamma \in \mathbf{F}_q$, $q = 2^e$, has trace 1, then the $\gamma$-twist of $E : y^2 + xy = x^3 + a_2 x^2 + a_6$ is $E' : y^2 + xy = x^3 + (a_2 x + \gamma) x^2 + a_6$, see [SZ03, p. 77].

Observe that a quadratic nonresidue becomes a quadratic residue in a quadratic extension of $\mathbf{F}_q$. Hence an elliptic curve $E$ and its twist $E'$ becomes isomorphic over $\mathbf{F}_{q^{2n}}$. However, the pair of twisted curves $E$ and $E'$ are not isomorphic over of $\mathbf{F}_{q^n}$, $n$ odd, even though these curves have the same $j$-invariant. In fact, for nonsupersingular curve $E(\mathbf{F}_q) = q + 1 \pm t$, for some $t \neq 0$. Thus their groups of points have distinct cardinalities.

**Theorem 23.** If the curve $E : f(x, y) = 0$ is the quadratic twist of the curve $E' : g(x, y) = 0$, then $\#E(\mathbf{F}_q) + \# E'(\mathbf{F}_q) = 2q + 2$.

Proof: Assume $p > 2$. For every point $P = (x, y)$ on the curve $E$, the twisted point $P' = (ux, u^2 y)$ is not on the curve $E'$. Conversely, for every point $P = (x, y)$ not on the curve $E$, the twisted point $P' = (ux, u^2 y)$ is on the curve $E'$. Therefore, $E(\mathbf{F}_q) = q + 1 - t$, and $E'(\mathbf{F}_q) = q + 1 + t$. ∎

The proof in characteristic $p = 2$ is slightly more intricate, see [EN99, p. 102].

## 14. Coordinates and Addition Formulas

The points of an elliptic curve can be represented in various forms: affine coordinates $P = (x, y)$, projective coordinates $P = (x, y, z)$, etc. The different choices of coordinates have specific and unique advantages and disadvantages such as simpler doubling, scalar multiplication, simpler formulas, faster implementations, etc., see [FT06, p. 280], [HM05] for more details.

**Group Law on $E(\mathbb{F}_q)$**
(1) $P + O = O + P$,              Identity,
(2) $P - P = -P + P$,             Additive inverse,
(3) $(P + Q) + R = P + (Q + R)$,  Associative,

The addition formulas can be written in various equivalent forms, and there are different formulas for the binary fields and nonbinary fields.

**Points Addition and Duplication on the Long Form Equation**
Let $P_1 = (x_1, y_1)$, $P_2 = (x_2, y_2)$, and $P_3 = (x_3, y_3)$ be points on the elliptic curve specified by the long form equation $E : y^2 + a_1 xy + a_3 y = x^3 + a_2 x^2 + a_4 x + a_6$. Then

(1) **Additive inverse.** $-P_1 = (x_1, -y_1 - a_2 x_1 - a_3)$
(2) **Additive identity.** If $P_1 = -P_2$, the first coordinates are equal $x_1 = x_2$, then $(x_3, y_3) = P_1 + P_2 = P_1 - P_1 = O$.
(3) **Duplication.** If $P_1 = P_2$, the first coordinates are equal $x_1 = x_2$, then

$$(x_3, y_3) = P_1 + P_2 = 2P_1 = (U^2 + a_1 U - a_2 - x_1 - x_2, -(U + a_1)x_3 - V - a_3), \text{ where}$$

$$U = \frac{3x_1^3 + 2a_2 x_1 + a_4 - a_1 y_1}{2y_1 + a_1 x_1 + a_3}, \quad V = \frac{-x_1^3 + a_4 x_1 + 2a_6 - a_3 y_1}{2y_1 + a_1 x_1 + a_3} = y_1 - x_1 U. \tag{21}$$

(4) **Additive formula.** $P_1 + P_2 = P_3$,      if $x_1 \neq x_2$ and $P_1 \neq -P_2$, then

$$(x_3, y_3) = (U^2 + a_1 U - a_2 - x_1 - x_2, -(U + a_1)x_3 - V - a_3), \text{ where}$$

$$U = \frac{y_2 - y_1}{x_2 - x_1}, \quad V = \frac{x_2 y_1 - x_1 y_2}{x_2 - x_1} = y_1 - x_1 U. \tag{22}$$

***Example 24.*** As an exercise, add and double the points $P_1 = (x_1, y_1)$, $P_2 = (x_2, y_2)$, on the elliptic curve $E : y^2 + xy + 2y = x^3 + 3x^2 + 4x + 5$ over the rational numbers and over a finite field.

**Points Addition and Duplication on Short Form Equation**
The short form equation has simpler duplication and addition formulas, these are obtained from the long form equation by setting $a_1 = a_2 = a_3 = 0$, $a_4 = a$, $a_6 = b$. Let $P_1 = (x_1, y_1)$, $P_2 = (x_2, y_2)$, and $P_3 = (x_3, y_3)$ be points on the elliptic curve specified by the short form equation $E : y^2 = x^3 + ax + b$. Then

(1) **Duplication.** If $P_1 = P_2$, the first coordinates are equal $x_1 = x_2$, then

$$(x_3, y_3) = 2P_1 = (U^2 - 2x_1, -Ux_3 - V), \text{ where}$$

$$U = \frac{3x_1^3 + a}{2y_1}, \quad V = \frac{-x_1^3 + a_4 x_1 + 2a_6}{2y_1} = y_1 - x_1 U. \quad (23)$$

(2) **Additive formula.** $P_1 + P_2 = P_3$,      if $x_1 \neq x_2$ and $P_1 \neq -P_2$, then

$(x_3, y_3) = (U^2 - x_1 - x_2, -Ux_3 - V)$, where

$$U = \frac{y_2 - y_1}{x_2 - x_1}, \quad V = \frac{x_2 y_1 - x_1 y_2}{x_2 - x_1} = y_1 - x_1 U, \quad (24)$$

**Fast Scalar Point Multiplication**
The top computational problem in the implementation of elliptic curve cryptography is the computations of scalar multiples $nP$. The choice of coordinates is a major factor in the overall complexity of the algorithm. Fast exponentiation methods in multiplicative groups, and its equivalent scalar multiplication in additive groups, are active areas of research. The basic prototype of a scalar point multiplication is given here.

*Binary Method of Scalar Point Multiplication*
The additive version of the classical *square-and-multiply* algorithm used to compute powers of elements in a multiplicative group is called the *double-and-add* algorithm. The *double-and-add* algorithm is employed to compute scalar multiple of points in the additive group $E(\mathbf{F}_q)$ in the following manner:

$$nP = 2a_1(\cdots + (2a_{k-2}(2a_k P + a_{k-1} P) + a_{k-3} P) + a_{k-3} P) + \cdots) + a_0 P, \quad (25)$$

where $n = a_k 2^k + \cdots + a_1 2 + a_0$, $a_i \in \{0, 1\}$, is the binary expansion of the integer $n$. This represents the sum $nP = P + P + \cdots + P$, with $P, nP \in E(\mathbf{F}_q)$.

There are various scalar point multiplication methods. Some of these methods improve on the basic *double-and-add* algorithm by replacing the binary expansion of the integer $n$ with some other integer expansion. For example, *non-adjacent form* integer representations, addition chain, etc. A non-adjacent form integer representation is a signed binary expansion $n = a_k 2^k + \cdots + a_1 2 + a_0$, $a_i \in \{-1, 0, 1\}$, such that no pair of consecutive $a_i$, $a_{i+1}$ are nonzeros, for example, $a_i a_{i+1} \neq 0$ does not hold for all $i$. Each integer has a unique NAF representation, with an expected $n/3$ nonzero digits $a_i$, see [ME97], [ST07]. Other advanced techniques are known too, see [SA00], [MR98].

## 15. *n*th Division Polynomials
Let $P = (x, y)$ be a point on the elliptic curve $E : y^2 = x^3 + ax + b$. The scalar multiple $[n]P$ of points is a rational function of the coordinates $x, y$. More precisely

$$[n]P = \left( \frac{\phi(x, y)}{\psi(x, y)^2}, \frac{\omega(x, y)}{\psi(x, y)^3} \right), \quad (26)$$

where $\phi(x, y)$, $\psi(x, y)$, $\omega(x, y)$ are polynomials. The points of order $n$ are the roots of the *nth division polynomial* $\psi(x, y)$. This is indicated by the equation $\psi(x, y) = 0$, so that multiple $[n]P = (\infty, \infty)$. The $n$th division polynomials $\psi(x, y)$ are important in points related algorithms.

The division polynomials are specific polynomial functions of two variables defined on elliptic curves. The recursive construction is as follows.

$\psi_{-1}(x, y) = -1,$ $\quad\psi_2(x, y) = 2y,$
$\psi_0(x, y) = 0,$ $\quad\psi_3(x, y) = 3x^4 + 6ax^2 + 12bx - a^2,$
$\psi_1(x, y) = 1,$ $\quad\psi_4(x, y) = 4y(x^6 + 5ax^4 + 20bx^3 - 5a^2x^2 - 4abx - 8b^2 - a^3),$

$$\psi_{2n}(x) = \psi_n[\psi_{n+2}\psi_{n-1}^2 - \psi_{n-2}\psi_{n+1}^2]/2y, \qquad \psi_{2n+1}(x) = \psi_{n+2}\psi_n^3 - \psi_{n-1}\psi_{n+1}^3,$$

(1) $\varphi_n(x) = x\psi_n^2(x) - \psi_{n+1}(x)\psi_{n-1}(x) \in \mathbb{Z}[a, b, x]$

(2) $4y\omega_n(x) = \psi_{n+2}(x)\psi_{n-1}^2(x) - \psi_{n-2}(x)\psi_{n+1}^2(x) \in \mathbb{Z}[a, b, x]$

The basic polynomials of the sequence are

$\varphi_n(x) = x^{n^2} +$ lower terms, and $\psi_n^2(x) = n^2 x^{n^2-1} +$ lower terms. Define the polynomial

$$f_n(x) = \begin{cases} \psi(x, y^2 = x^3 + ax + b) & \text{if } n = \text{odd}, \\ \psi(x, y)/y & \text{if } n = \text{even}, \end{cases} \qquad \deg f_n(x) = \begin{cases} (n^2 - 1)/2 & \text{if } n = \text{odd}, \\ (n^2 - 4)/2 & \text{if } n = \text{even}. \end{cases} \tag{27}$$

**Lemma 25.** Let $n > 0$, then $f_n(x) = 0 \Leftrightarrow nP = O$.

## 16. Some Public Keys Protocols
The elliptic curve cryptography standards are specified in NIST FIPS 186-2, IEEE1363, SECQ and other documents issued by other standards organizations. A few important criteria are described here.

**Selection of Elliptic Curves and Finite fields**
The elliptic curves used in cryptographic applications should be chosen at random to avoid weak curves. In addition, the background finite fields should be selected to meet strong structural properties. Cryptographic grade elliptic curves and finite fields should have the following properties:

1. The order $\#E(\mathbf{F}_q)$ of the group $E(\mathbf{F}_q)$ is divisible by a large prime.
2. A large subgroup of $E(\mathbf{F}_q)$ has large embedding degree $k > (\log q)^3$ in the extension $\mathbf{F}_{q^k}$.
3. The order of the group $\#E(\mathbf{F}_q) \neq q$, and $\#E(\mathbf{F}_q) = rn$, $r$ small and $n$ a large prime.
4. The nonweak finite field $\mathbf{F}_q$, for example, $q = p^n$, where $n$ is not smooth, (have large prime factors) etc.
5. Other criteria.

Weak curves have special structural properties which allow easier computations of the discrete logarithms on their groups of points. Similarly, weak finite fields have richer structures which allow easier computations of the discrete logarithms. The best known classes of weak curves and finite fields are the following.

1. The order of the group $\#E(\mathbf{F}_q)$ is divisible by small primes only.
2. Any subgroup of $E(\mathbf{F}_q)$ has small embedding degree $k < (\log q)^3$ in the extension $\mathbf{F}_{q^k}$.
3. Anomalous curve of group order $\#E(\mathbf{F}_q) = q$.
4. The weak finite field $\mathbf{F}_q$, for example, $q = p^n$, where $n$ is smooth, (have only small prime factors) etc.
5. Other criteria

Random curves generation and verification algorithms are readily available in the literature, see [NT99, Appendices 4 and 5], [SECQ] etc. The weak finite fields are discussed in [MT04].

**Text Encoding and Decoding,**
A plaintext encoding algorithm converts plaintext in the interval [0, M] to plaintext as a point in the group of points $E(\mathbf{F}_q)$.
*Plaintext Encode Algorithm* I.
Input: $m \in [0, M^2]$ such that $q > M(K + 1)$, and $E : y^2 = x^3 + ax + b$.
Output: A point $P = (x_0, y_0) = (mK + i_0, y_0)$ on the elliptic curve.
1. For $i < K$, While $(mK + i)^3 + a(mK + i) + b \ne y^2$.
2. Return the point $P = (mK + i_0, y_0)$.
This algorithm encode a plaintext value in the range 0 to M with a negligible probability of failure, namely, $2^{-K}$.
The decoding process is simply

$[\, x_0/K \,] = [(mK + i_0)/K] = m$, which is the largest integer after division by K.

**FIPS186-2 Standard**
The FIPS186-2 Standard lists 15 elliptic curves of various level of security.
(i) There are 5 random nonbinary ordinary elliptic curves. The parameters are specified as 6-tuples $(p, a = -3, b, h, n, P)$. Here $p$ is a prime, $b$ is a semi random integer modulo $p$, the group order $hn = \#E(\mathbf{F}_q)$ is prime, and $P$ is a point of order $n$. The equation of the curve is $E : y^2 = x^3 - 3x + b$ over the finite fields $\mathbf{F}_p$, where the primes are

$p_{192} = 2^{192} - 2^{64} - 1,$
$p_{256} = 2^{256} - 2^{224} + 2^{192} + 2^{96} - 1,$
$p_{521} = 2^{521} - 1.$

$p_{224} = 2^{224} - 2^{96} + 1,$
$p_{384} = 2^{384} - 2^{128} - 2^{96} + 2^{32} - 1,$

These elliptic curves are semirandom since one of the coefficients is not random. The pre-assignment of the coefficient $a = -3$ (or more generally $a = -3c^2$, $c \ne 0$) reduces the number of arithmetic operations used to compute $2P$ in projective coordinates by about 50%, see [HM04, p.90], [FT06, p. 282].

(ii) There are 5 random binary ordinary elliptic curves $E : y^2 + xy = x^3 + ax + b$ over $\mathbf{F}_{2^n}$, with $a$, $b = 0, 1$. The group order is $\#E(\mathbf{F}_q) = hr$, $r$ prime, and $h = 2, 4$.

(iii) There are 5 random binary anomalous curves $E : y^2 + y = x^3 + ax + b$ over $\mathbf{F}_{2^n}$, with $a$, $b = 0, 1$. The group order is $\#E(\mathbf{F}_q) = hr$, $r$ prime, and $h = 1$.

In both cases the binary finite fields use polynomial representations $\mathbf{F}_{2^n} \cong \mathbf{F}_2[x]/(f(x))$, where

$f(x) = x^{163} + x^7 + x^6 + x^3 + 1,$
$f(x) = x^{283} + x^{12} + x^7 + x^5 + 1,$
$f(x) = x^{571} + x^{10} + x^5 + x^2 + 1.$

$f(x) = x^{233} + x^{74} + 1,$
$f(x) = x^{409} + x^{87} + 1,$

## 17. Counting the Number of Points on Curves
The computations of the group of points $E(\mathbf{F}_q)$ and the order $\#E(\mathbf{F}_q)$ are two different problems. The former seeks a list of all points, which for most curves of interest in applications this is not a practical task. In contrast, the later seeks the number of points in the set, this is a more feasible task.

For each $x \ne 0$, the polynomial $f(x) = x^3 + ax^2 + bx + c$ can either be a zero, a square or a nonsquare in the finite field $\mathbf{F}_q$. Accordingly, the equation $y^2 = f(x)$ has either 1 solution, 2 solutions or no solution. As a consequence, the set of rational points $E(\mathbf{F}_q) = \{ (x, y) \in \mathbf{F}_q \times \mathbf{F}_q : y^2 = f(x) \}$ contains at most $2q$ points, this is a trivial estimate. The order of the group of points $E(\mathbf{F}_q)$ is an integer $N$ in the Hasse interval

$$[q + 1 - 2q^{1/2}, q + 1 + 2q^{1/2}]. \tag{28}$$

Each point $P \in E(\mathbf{F}_q)$ satisfies the order equation $NP = O$, where $O$ is the identity point.

**Theorem 26.** (Hasse 1934). Let $E : y^2 = f(x, y)$ be an elliptic curve over $\mathbf{F}_q$. Then
(1) The number of points $N_n$ on a curve of genus $g \geq 1$ satisfies the inequality $| q^n + 1 - N_n | \leq 2gq^{1/2}$.
(2) The kernel of the endomorphism $\sigma - 1$ is $\ker(\sigma - 1) = E(\mathbf{F}_q)$, where $\sigma$ is the Frobenious map.

The original proof of this result was given a single line as opposed to the lengthy proof usually given, see [LG05, p. 619]. A proof based on the estimate of the corresponding cubic exponential sum is also straight forward and short.

**Definition 27.** A curve $C : f(x, y) = 0$ over $\mathbf{F}_q$ of genus $g > 0$ is called *maximal curve* (minimal) if it attains the maximal (minimal) number of points $\#C(\mathbf{F}_q) = q^n + 1 \pm 2gq^{1/2}$.

Clearly, maximal (minimal) curves exist only over quadratic extensions, that is, for $q = p^{2n}$, $n \geq 1$.

The best known methods of determining the order $\#E(\mathbf{F}_q)$ of the group of points $E(\mathbf{F}_q)$ are the following:
1. Brute force calculation,
2. Baby step giant step method,
3. Recursive counting, and manipulation of the zeta function,
4. Closed form evaluations,
5. Advanced techniques.

(1) **Brute Force Calculations.** The brute force enumeration of the group of points $E(\mathbf{F}_q) = \{ (x, y) : y^2 = f(x)$ with $x, y \in \mathbf{F}_q \}$ can be accomplished by listing all the points $(x_1, y_1), (x_2, y_2), \ldots, (x_n, y_n)$ such that $y^2 - f(x) = 0$. The naïve search algorithm runs in $O(q)$ finite fields operations. This is practical for very small prime powers $q$, depending on current computing technology.

A more efficient brute force enumeration algorithm that runs in about $O(q^{1/2}(\log q)^4)$ finite fields operations emerges from Hasse's theorem:

*Group Order Computation Algorithm* I
1. Input: A nonsingular curve $E : y^2 = f(x)$ over $\mathbf{F}_q$.
2. Select a random point $P \in E(\mathbf{F}_q)$ on the curve.
3. If $mP = O$ for some $m \leq 4q^{1/2}$, then repeat step 2.
4. Compute the integer $n$ such that $nP = O$, and $q + 1 - 2\sqrt{q} \leq n \leq q + 1 + 2\sqrt{q}$.
5. Output: The order $n = \# E(\mathbf{F}_q)$.

Step 3 in the algorithm eliminates points of small orders $n \leq 4q^{1/2}$, this prevents the existence of two or more integers $n_0, n_1, \ldots, n_d \geq q + 1 - 2q^{1/2}$ such that $n_0 P_0 = n_1 P_1 = \cdots n_d P_d = O$. For example, if $2P = O$, then every integer $n = 2k \geq q + 1 - 2q^{1/2}$ satisfies $nP = O$, if $3P = O$, then every integer $n = 3k \geq q + 1 - 2q^{1/2}$ satisfies $nP = O$, and so on. But requiring that $nP = O$ if and only if $n > 4q^{1/2}$, implies that the order computed by the algorithm is unique.

**Theorem 28.** (Mestre)  Let $p > 229$. Every integer $N \in [q + 1 - 2\sqrt{q}, q + 1 + 2\sqrt{q}]$ is the unique order of either an elliptic curve $E$ or its quadratic twist $E'$ over $\mathbf{F}_q$.

**(2) Baby Step Giant Step Method**
The Baby step giant step is a general technique for computing certain parameters of discrete structures in certain ranges. The parameters can be the orders of elements in cyclic groups, the discrete logarithms, and so on. The space and time complexities vary depending on the version implemented, but each one requires approximately $O(N^{1/2+\varepsilon})$ space allocations and group operations on a group $G$ of order $\#G = N$.

The application of interest here is the calculation of the trace of an elliptic curve, which is equivalent to the calculation of the order of its group of points.

Put $Q = (q + 1 + [2q^{1/2}])P$, where $[\,x\,]$ denotes the largest integer function. The basic idea in the determination of the group order of an elliptic curve is the observation that

$$Q = (q + 1 + [2q^{1/2}])P = (q + 1 - t + t + [2q^{1/2}])P = (t + [2q^{1/2}])P = nP, \tag{29}$$

for some $0 \le n < 4q^{1/2}$ or equivalently $0 \le |t| < 2q^{1/2}$. This procedure requires a point $P$ of order $\mathrm{ord}(P) > 4q^{1/2}$ in order to have a unique solution. The algorithm finds an integer $n \in [0, 4q^{1/2}]$ such that $Q = nP$, with $n = am + b$, $0 \le a, b < 2q^{1/2}$, $m = 2q^{1/4}$. After the integer $n$ is found, it returns the order $N = \#E(\mathbf{F}_q) = q + 1 - t$ with $t = am + b - 2q^{1/2}$.

Since $|t| < 2q^{1/2}$, the sorting takes approximately $O(N^{1/4+\varepsilon})$ operations and the same for finding a match. Hence, the running time of the algorithm is approximately $O(N^{1/4+\varepsilon})$ operations in $E(\mathbf{F}_q)$.

*BSGS Algorithm*
Input: Input: A nonsingular curve $E : y^2 = f(x)$ over $\mathbf{F}_q$.
Output: The order $N = \#E(\mathbf{F}_q) = q + 1 - t$.
1. Select a random point $P \in E(\mathbf{F}_q)$ on the curve of order $\mathrm{ord}(P) > 4q^{1/2}$.
2. Put $Q = (q + 1 + [2q^{1/2}])P$,
3. Baby steps. Store the points $P, 2P, 3P, \ldots, bP, \ldots, kP$, with $k = [2q^{1/4}]$.
4. Giant steps. Store the points $Q - akP$, for $a = 1, 2, \ldots, k$.
5. Sort the Giant steps list to find a match $Q - akP = bP$, in the Baby steps list.
6. Return $t = ak + b$, and $\#E(\mathbf{F}_q) = q + 1 - t$.

**(3) Recursive Counting**
The recursive computation of the number of rational points can be accomplished in about three methods.
(i) Quadratic Recursive Formulas.
(ii) Repeated differentiation of the zeta function of the curve.
(iii) Modular Form Recursive Formulas, see Example 77.

An algorithm based on quadratic recurrent sequence is described here. This algorithm springs from the decomposition formula $\#E(\mathbf{F}_q) = (\pi - 1)(\bar{\pi} - 1) = q + 1 - \pi - \bar{\pi}$. This is based on the fact that if $q = \pi\bar{\pi}$ in $\mathrm{End}(E)$, then $\pi = t + u\sqrt{d}$ is an eigenvalue of the Frobenious map. Specifically, $(x, y) \to (x^q, y^q) = [\pi](x, y)$, and the theory of recurrent sequences.

The theory of quadratic recurrent sequences states that for a given $f(x) = x^2 - Px + Q$ of discriminant $D = P^2 - 4Q$, with a pair of conjugate roots $\alpha, \beta$, there is a natural pair of second order recurrent sequences

$$U_n(P,Q) = \frac{\alpha^n + \beta^n}{\alpha - \beta}, \quad \text{and} \quad V_n(P,Q) = \alpha^n + \beta^n. \tag{30}$$

The later sequence satisfies the recurrent relation

(i) $V_n = V_1 V_{n-1} - QV_{n-2}$, $n \geq 2$, and the replication formulas

(ii) $V_{2n} = V_n^2 - 2Q^n$, (iii) $V_{3n} = V_n(V_n^2 - 3Q^n)$, (iv) $V_{m+n} = V_m(V_n - QV_{m-n})$,

see the literature for more details. These identities are used below to compute the integer $N_{n+1} = \#E(\mathbf{F}_{q^{n+1}})$, given $N_n = \#E(\mathbf{F}_{q^n})$ recursively. The replication formulas (ii) to (iv) are useful for computing the orders of nonconsecutive groups of points such as $N_n = \#E(\mathbf{F}_{q^n})$ and $N_{3n} \#E(\mathbf{F}_{q^{3n}})$.

***Theorem 29.*** Let $\#E(\mathbf{F}_q) = q + 1 - t_1$, and let $n \geq 2$. Then the order of the group of points $\#E(\mathbf{F}_{q^n}) = q^n + 1 - V_n$, where $V_n = t_1 V_{n-1} - qV_{n-2}$, $V_0 = 2$, $V_1 = t_1$.

Proof: Since $\#E(\mathbf{F}_{q^n}) = q^n + 1 - \alpha^n - \overline{\alpha}^n$, and $\alpha, \beta = \overline{\alpha}$ are the roots of $L(T) = T^2 - t_1 T + q$, the term $V_n = \alpha^n + \overline{\alpha}^n$ satisfies the specified quadratic recurrent sequence. ∎

***Definition 30.*** The infinite sequence $\{ N_n : n \geq 1 \}$ is a *divisibility sequence*. More precisely, if $d \mid n$, then $N_d \mid N_n$.

In terms of subgroups of $\mathbf{F}_{q^n}$, this states that, if $d$ is a divisor of $n$, then $\#E(\mathbf{F}_{q^d})$ is a subgroup of $E(\mathbf{F}_{q^n})$.

***Example 31.*** The elliptic curve $E : y^2 + y = x^3$ and its group of points $E(\mathbf{F}_2) = \{ O, (0, 0), (0, 1) \}$, here $\#E(\mathbf{F}_2) = 3$. This is a Supersingular curve since its characteristic polynomial $f_\pi(x) = x^2 + 2$. On the next larger field $\mathbf{F}_4$ it has $\#E(\mathbf{F}_4) = 2^2 + 1 - V_2 = 9$ points since $V_2 = (V_1)^2 - qV_0, = -4$, the trace $t_1 = 0$. The numbers $N_n = \#E(\mathbf{F}_{2^n})$ are tabulated for $n < 16$.

| N | $V_n$ | $N_n$ | n | $V_n$ | $N_n$ | n | $V_n$ | $N_n$ |
|---|---|---|---|---|---|---|---|---|
| 1 | 0 | 3 | 6 | −16 | $3^4$ | 11 | 0 | 3·683 |
| 2 | −4 | $3^2$ | 7 | 0 | 3·43 | 12 | 128 | $3^4 \cdot 7^2$ |
| 3 | 0 | $3^2$ | 8 | 32 | $3^2 \cdot 5^2$ | 13 | 0 | 3·2731 |
| 4 | 8 | $3^2$ | 9 | 0 | $3^2 \cdot 19$ | 14 | −512 | $3^2 \cdot 43^2$ |
| 5 | 0 | 3·11 | 10 | −64 | $3^2 \cdot 11^2$ | 15 | 0 | $3^2 \cdot 11 \cdot 331$ |

Some of these extensions with $n > 160$ which have large prime factors $p \mid N_n$, are suitable for cryptographic applications.

A recursion formula for curves of genus $g > 1$, which generalizes this technique, can be derived using the Newton-Girard formula.

***Example 32.*** This example demonstrates another way of computing $\#E(\mathbf{F}_q)$ recursively. As an illustration, take the elliptic curve $E : y^2 = x^3 + x$ over $\mathbf{F}_5$ and its group of points $E(\mathbf{F}_5) = \{ O, (0, 0), (\pm 2, 0) \}$, here $\#E(\mathbf{F}_5) = 4$. Since $\text{End}(E) = \mathbb{Z}[i]$, and the prime splits as $p = 5 = (1 + i2)(1 - i2)$, the (eigenvalue) Frobenious map is $\pi = 1 + i2$, and $(\pi - 1)(\overline{\pi} - 1) = q + 1 - \pi - \overline{\pi} = 4$. Over the next larger field $\mathbf{F}_{25}$ it has $\#E(\mathbf{F}_{5^2}) = 5^2 + 1 - \pi^2 - \overline{\pi}^2 = 32$ points, and so on. The numbers $N_n = 5^n + 1 - \pi^n - \overline{\pi}^n$ are tabulated for $n < 11$.

| N | $\pi^n$ | $N_n$ | n | $\pi^n$ | $N_n$ |
|---|---|---|---|---|---|
| 1 | $1 + 2i$ | $2^2$ | 6 | $117 + 44i$ | $2^5 \cdot 13 \cdot 37$ |
| 2 | $-3 + 4i$ | $2^5$ | 7 | $29 + 278i$ | $2^2 \cdot 29 \cdot 673$ |
| 3 | $-11 + 2i$ | $2^2 \cdot 37$ | 8 | $-527 + 336i$ | $2^9 \cdot 3^2 \cdot 5 \cdot 7 \cdot 17$ |
| 4 | $-7 - 24i$ | $2^7 \cdot 5$ | 9 | $-1199 - 718i$ | $2^2 \cdot 37 \cdot 73 \cdot 181$ |
| 5 | $41 - 38i$ | $2^2 \cdot 761$ | 10 | $237 - 3116i$ | $2^5 \cdot 401 \cdot 761$ |

## 18. Closed Form Techniques

The basic derivation of the exponential sum involved in the calculation of the order of the group $E(\mathbf{F}_q)$ is as follows: Let the quadratic symbol be defined by

$$\left(\frac{x}{\mathbf{F}_q}\right) = \begin{cases} 1 & \text{if } x \text{ is a square,} \\ 0 & \text{if } x \text{ is zero in } \mathbf{F}_q, \\ -1 & \text{if } x \text{ is not a square.} \end{cases} \tag{31}$$

Each $x \in \mathbf{F}_q$ contributes $1 + \left(\frac{f(x)}{\mathbf{F}_q}\right) = 1, 2$ or $0$ points ( $(x_0, 0)$, or $(x_0, \pm y_0)$, no point) to the cardinality $\#E(\mathbf{F}_q)$ of the group $E(\mathbf{F}_q)$. This occurs whenever the value of $f(x)$ at $x_0$ is a zero, a quadratic residue or a quadratic nonresidue respectively. Now summing over all the $x \in \mathbf{F}_q$ gives the formula

$$\#E(\mathbf{F}_q) = 1 + \sum_{x \in \mathbf{F}_q} \left(1 + \left(\frac{f(x)}{\mathbf{F}_q}\right)\right) = q + 1 + \sum_{x \in \mathbf{F}_q} \left(\frac{f(x)}{\mathbf{F}_q}\right), \tag{32}$$

for the order of $E(\mathbf{F}_q)$ of any elliptic curve $E : y^2 = f(x)$. The count includes the point at infinity.

In practice this technique is limited to special curves for which the exponential sum can be evaluated easily, and to small finite fields. There are various cubic exponential sums that have exact evaluations. Some of these cases are the 9 CM elliptic curves over the rational numbers $\mathbb{Q}$. Under these conditions the prime $p$ splits as $p = \pi\bar{\pi}$, $\pi = (u + v\sqrt{d})/2$ in $\mathbb{Q}(\sqrt{d})$, $d = -1, -2, -3, -7, -11, -19, -43, -67, -163$, then $L(T) = T^2 - a_p T + p$, where $a_p = -\text{Tr}(\pi) \neq 0$. This event occurs whenever the quadratic symbol $\left(\frac{d}{p}\right) = 1$. In addition, $E(\mathbf{F}_p) = p + 1 + \sum_{x \in \mathbf{F}_p} \left(\frac{f(x)}{p}\right) = p + 1 - a_p$. Thus the exponential sum $S(f(x)) = \sum_{x \in \mathbf{F}_p} \left(\frac{f(x)}{p}\right) = -a_p$. Otherwise, $\left(\frac{d}{p}\right) = -1$ and the prime $p$ is inert in $\mathbb{Q}(\sqrt{d})$. Then $L(T) = T^2 + p$, and $a_p = -\text{Tr}(\pi) = 0$.

The calculations of the orders for two of these curves are well known, they are given here.

**Proposition 33.** Let $q = 3m + 2$, and let $E : y^2 = x^3 + b$. Then $\#E(\mathbf{F}_q) = q + 1$.

Proof: The correspondence $x \to x^3$ rearranges all the elements of $\mathbf{F}_q$ in a one-to-one manner. Likewise, $x \to x^3 + a$, $a \in \mathbf{F}_q$. Thus

$$E(\mathbf{F}_q) = 1 + \sum_{x \in \mathbf{F}_q}\left(1 + \left(\frac{x^3+b}{\mathbf{F}_q}\right)\right) = q + 1 + \sum_{x \in \mathbf{F}_q}\left(\frac{x^3+b}{\mathbf{F}_q}\right) = q + 1 + \sum_{x \in \mathbf{F}_q}\left(\frac{x}{\mathbf{F}_q}\right) = q+1. \tag{33}$$

The steps taken there use the fact that there are equal amounts of quadratic residues and nonquadratic residues in $\mathbf{F}_q$, and the value of the quadratic symbol $(x \mid \mathbf{F}_q) = 1$ or $-1$ whenever $x$ is a quadratic or nonquadratic residue respectively. ∎

The endomorphism ring of the elliptic curve $E : y^2 = x^3 + b$ is $\text{End}(E) = \mathbb{Z}[\omega]$, $\omega^3 = 1$, but $\omega \neq 1$. The scalar multiple $[\omega]P = (\omega x, y)$ for a given point $P = (x, y)$. For an arbitrary prime $\pi \in \mathbb{Z}[\omega]$ not divisible by 6, the order of the group is $\#E(\mathbf{F}_\pi) = N(\pi) + 1 - \text{Tr}(\pi)$. These observations partially prove the following, the complete details are in the literature.

***Proposition 34.*** Let $E : y^2 = x^3 + b$, and $p = u^2 + 3v^2$, $u \equiv 2 \bmod 3$. Then

(i) $\#E(\mathbf{F}_p) = \begin{cases} p+1+2u & \text{if } p \equiv 1 \bmod 4, \\ p+1-2u & \text{if } p \equiv 3 \bmod 4. \end{cases}$

(ii) $\#E(\mathbf{F}_{p^2}) = \begin{cases} (p+1)^2 & \text{if } p \equiv 1 \bmod 4, \\ (p-1)^2 & \text{if } p \equiv 3 \bmod 4. \end{cases}$

***Proposition 35.*** Let $E : y^2 = x^3 + ax$, and $p \equiv 3 \bmod 4$, $\gcd(a, p) = 1$. Then

$$\#E(\mathbf{F}_{p^n}) = \begin{cases} p^n + 1 & \text{if } n \equiv 1 \bmod 2, \\ p^n + (-1)^{n/2} 2p^{n/2} & \text{if } n \equiv 0 \bmod 2. \end{cases}$$

Proof : For $n \equiv 1 \bmod 2$, use the fact that either $x$ or $-x$ is a quadratic residue in characteristic $p \equiv 3 \bmod 4$, everything else is similar to the previous one. For $n \equiv 0 \bmod 2$, use a root $\alpha$ of the characteristic polynomial $L(T) = T^2 + p$ to compute $N_n = p^n + 1 - \alpha^n - \overline{\alpha}^n$. ∎

This analysis given above appears in [SZ03, p.68] and similar sources. Observe that exactly the same analysis applies to any one-to-one cubic polynomial $f(x) = x^3 + ax^2 + bx + c$ over $\mathbf{F}_q$ besides $f(x) = x^3 + ax$ or $x^3 + b$. A related result is also included here.

***Proposition 36.*** Let $E : y^2 = x^3 + ax$, and $p = u^2 + v^2$, $u \equiv -1 \bmod 4$. Then

$$\#E(\mathbf{F}_p) = \begin{cases} p+1-2u & \text{if } a^{(p-1)/4} \equiv 1 \bmod p, \\ p+1+2u & \text{if } a^{(p-1)/4} \equiv -1 \bmod p, \\ p+1-2v & \text{if } a^{(p-1)/4} \equiv -u/v \bmod p, \\ p+1+2v & \text{if } a^{(p-1)/4} \equiv u/v \bmod p. \end{cases}$$

This result attributed to Gauss has several proofs, based on the evaluation of the cubic exponential sum.

## 19. Hasse Invariant
Given a prime $p = 2n + 1$, the Hasse polynomial (or Hasse invariant) is defined by

$$H_p(x) = (-1)^n \sum_{k=0}^{n} \binom{n}{k}^2 x^k. \tag{8}$$

*Examples*
$H_3(x) = -(x + 1) \in \mathbf{F}_3[x]$, For $\beta = -1$ the *j*-invariant is $j(E) = 0$, which is supersingular over $\mathbf{F}_3$.

$H_5(x) = x^2 + 2x + 1 = (x + 1)^2 \in \mathbf{F}_5[x]$. For $\beta = -1$ the *j*-invariant is $j(E) = 2$, which is supersingular over $\mathbf{F}_5$.

$H_7(x) = x^3 + 2x^2 + 2x + 1 = (x + 1)(x - 2)(x + 3) \in \mathbf{F}_7[x]$. For $\beta = -1, 2, 3$, the *j*-invariant is $j(E) = -1, 0, 1$, which is supersingular over $\mathbf{F}_7$.
$H_{11}(x) = x^5 + 3x^4 + x^3 + x^2 + 3x + 1 = (x^2 - x + 1)(x + 1)(x - 2)(x + 5) \in \mathbf{F}_{11}[x]$. For $\beta = -1, 2, -5$, the *j*-invariant is $j(E) = 0, 1, -1$ which is supersingular over $\mathbf{F}_{11}$.

$H_{13}(x) = x^6 - 3x^5 - 4x^4 + 3x^3 - 4x^2 - 3x + 1 = (x^6 + ax + b)(x - 4)(x - 10) \in \mathbf{F}_{13}[x]$. For $\beta = -1, 2, -5$, the *j*-invariant is $j(E) = 0, 1$, which is supersingular over $\mathbf{F}_{13}$.

The Hasse invariant is used to evaluate the character sum connected with the polynomial $f_\beta(x) = x(x - 1)(x - \beta)$, $\beta \neq 0, 1$. Specifically,

$$S(f_\beta) = \sum_{x=0}^{p-1} \left( \frac{x(x-1)(x-\beta)}{p} \right) \equiv (-1)^m \sum_{k=0}^{m} \binom{m}{k}^2 \beta^k \bmod p. \tag{8}$$

**Theorem 37.** (Manin 1965) Let $E: y^2 = f_\beta(x) = x(x - 1)(x - \beta)$, $0, 1 \neq \beta$, be an elliptic curve over the rational numbers $\mathbb{Q}$ with good reduction at a prime $p$, and let $a_p = -S(f_\beta)$ be the trace of Frobenius at $p$, then $H(\beta) \equiv a_p \bmod p$. In particular,

$$\#E(\mathbf{F}_p) = \#\{ (x, y) : y^2 - f_\beta(x) = 0, \text{ and } x, y \in \mathbf{F}_p \} = p + 1 + (H(\beta) \bmod p).$$

More generally, for an elliptic curve $E : y^2 = f_\beta(x)$, the cardinality of its group of points satisfies the congruence $\#E(\mathbf{F}_q) \equiv 1 - H(\beta)^{(q-1)/(p-1)} \bmod p$.

For $p > 2$, a supersingular Legendre form is characterized by the condition $H(\beta) \equiv 0 \bmod p$, see [HR87, p. 252].

**(3) Advanced Techniques**
Advanced techniques are all the new algorithms developed since the late 1980's. These new methods include Schoof Algorithm, and related techniques, Satoh algorithm, Kedlaya algorithm, and Mestre AGM algorithm. These algorithms are quite complex, see [FR05] for comprehensive introductions.

## 20. The Structure Of The Group $E(\mathbf{F}_q)$
The important problems in points counting are the order of the group and the structure of the group.

**Theorem 38.** (Cassels 1966). The group of points $E(\mathbf{F}_q)$ of a curve $E : y^2 = x^3 + ax^2 + bx + c$ is isomorphic to the group $\mathbb{Z}_m \times \mathbb{Z}_n$, with $m \mid n$, and $m \mid q - 1$.

**Proposition 39.** Let $E : y^2 = x^3 + ax^2 + bx + c$ be a nonsingular curve. Then the group $E(\mathbf{F}_q)$ is odd if and only if $f(x) = x^3 + ax^2 + bx + c$ is irreducible over $\mathbf{F}_q$.

Proof: It is known that $E(\mathbf{F}_q) \cong \mathbb{Z}_m \times \mathbb{Z}_n$, with $m \mid n$. Further, a zero $(x, 0)$ is a point of order 2, so the parity of $\#E(\mathbf{F}_q)$ is given by

$$\#E(\mathbf{F}_q) \equiv \begin{cases} 0 \bmod 4 & \text{if } f(x) \text{ has 3 zeros in } \mathbf{F}_q, \\ 2 \bmod 4 & \text{if } f(x) \text{ has 1 zeros in } \mathbf{F}_q, \\ 1, 3 \bmod 4 & \text{if } f(x) \text{ has no zeros in } \mathbf{F}_q. \end{cases} \quad (34)$$

Thus, if $f(x)$ is irreducible the parity is $\#E(\mathbf{F}_q)$ is odd and has no points of order 2. This implies that $m = 1$. ∎

The modulo $p$ reduction map $E(\mathbb{Q})_{\text{tor}} \to E(\mathbf{F}_q)$, which is group morphism, see [HR87] can be used to determine whether the group order $\#E(\mathbf{F}_q)$ is divisible by some $d \leq 16$. This limit $d = 16$ arises from the possible order of the torsion group $E(\mathbb{Q})_{\text{tor}}$. This is known to be $d \leq 16$.

**Theorem 40.** (Masur 1988) The torsion group of an elliptic curve is either $E(\mathbb{Q})_{\text{tor}} \cong \mathbb{Z}_n$ with $n = 1, 2, \ldots, 10, 12$ or $E(\mathbb{Q})_{\text{tor}} \cong \mathbb{Z}_2 \times \mathbb{Z}_n$ with $n = 2, 4, 6, 8$.

The order of the torsion group $E(\mathbb{Q})_{\text{tor}}$ over the rational number is not divisible by 11, 13 nor 14, however $\#E(\mathbf{F}_q) \equiv 0 \bmod 11, 13$ or 14 do occur.

For super singular elliptic curves, the torsion group is $E(\mathbb{Q})_{\text{tor}} \cong \mathbb{Z}_n$ with $n = 1, 2, 3, 4$ or 6, see [MY91] for an easy proof.

**Proposition 41.** (i) The group $E(\mathbf{F}_q)$ of points of $E : y^2 = x^3 + ax$ contains a subgroup isomorphic to $\mathbb{Z}_2$ or $\mathbb{Z}_4$ or $\mathbb{Z}_2 \times \mathbb{Z}_2$ if $a = 4$ or $a = -4$ or $a = -u^2$, respectively.

(ii) The group $E(\mathbf{F}_q)$ of points of $E : y^2 = x^3 + b$ contains a subgroup isomorphic to $\mathbb{Z}_1$ or $\mathbb{Z}_2$ or $\mathbb{Z}_3$ or $\mathbb{Z}_6$.

A discussion appears in [HR87, p.34]. The group of points of an elliptic curve that has three rational roots cannot be cyclic because $\mathbb{Z}_2 \times \mathbb{Z}$ is contained in $E(\mathbf{F}_q)$.

**Proposition 42.** Let $E : y^2 = x^3 + ax + b$ be a curve over $\mathbf{F}_q$, char $p > 3$. Then
(1) If $a, b \neq 0$, then for each $x \neq 0$ there are 0 or 2 points $(x, \pm y)$ on the curve.
(2) If $a \neq 0$ is a square, $b = 0$, and $p = 4n + 1$, then for each $x \neq 0$ there are 0 or 4 points $(x, \pm y)$ and $(-x, \pm cdy)$ on the curve, where $a = c^2$, and $-1 = d^2$. In particular $E(\mathbf{F}_q) \equiv 0 \bmod 4$.

Proof of 2: In characteristic $p = 4n + 1$, the equation $-1 = d^2$ is solvable. Accordingly $y^2 = x^3 + ax + b = x(x^2 + a) = x(x^2 - c^2d^2) = x(x - cd)(x + cd)$ has 3 roots in $\mathbf{F}_q$. Further, for $x \neq 0$ such that $(x, \pm y)$ is a points on the curve, and so is $(-x, \pm dy)$. This implies that $E(\mathbf{F}_q)$ contains a subgroup $G = \{O, (0, 0), (\pm cd, 0)\}$ of order 4. ∎

**Theorem 43.** (Waterhouse 1969). Let $p$ be a prime and let $E$ be an elliptic curve over $\mathbf{F}_q$, $q = p^n$, $n \geq 1$. Then the following hold.
i. $E$ is not supersingular and $\gcd(t_n, p) = 1$.
ii. $E$ is supersingular and $\gcd(t_n, p) \neq 1$.
iii. If $n$ is even and $p = 3m + 1$, then $t_n = \pm 2p^{n/2}$.
iv. If $n$ is even and $p = 3m + 2$, then $t_n = \pm p^{n/2}$.
v. If $n$ is odd and $p = 2$ or 3, then $t_n = \pm p^{(n+1)/2}$.
vi. If $n$ is odd and $p = 4a + 1$, or $n$ is even and $p = 4m + 3$, then $t_n = 0$.

***Theorem 44.*** (Schoof 1985). The possible structures of the group of points of supersingular elliptic curves are as follow.

1. If $t_n = 0$ and $n$ is odd, or $q \equiv 3 \bmod 4$, then $E(\mathbf{F}_q) \cong \mathbb{Z}_2 \times \mathbb{Z}_{(q+1)/2}$, or cyclic, otherwise cyclic.
2. If $t_n = \pm p^{n/2}$, or $q \equiv 2 \bmod 3$, then $E(\mathbf{F}_q) \cong \mathbb{Z}_N$ is cyclic.
3. If $t_n = \pm 2^{n/2+1}$, then $E(\mathbf{F}_q) \cong \mathbb{Z}_m \times \mathbb{Z}_m$, where $m = q^{1/2} \pm 1$.
4. If $t_n = \pm 2^{(n+1)/2}$, then $E(\mathbf{F}_q) \cong \mathbb{Z}_N$ is cyclic.
5. If $t_n = \pm 3^{(n+1)/2}$, then $E(\mathbf{F}_q) \cong \mathbb{Z}_N$ is cyclic.

These two results provide a broad classification of the structures of all CM elliptic curves. The complete structures of nonCM elliptic curves were later completed independently by two different authors.

***Theorem 45.*** ([RK87], [VL88]). If $t$ is an integer with $|t| \le 2q^{1/2}$, and $\gcd(t, q) = 1$, the possible groups that an elliptic curve over $\mathbf{F}_q$ with $N = q + 1 - t$ can be are

$$E(\mathbf{F}_q) \cong \mathbb{Z}_{p^{v_p(N)}} \oplus \mathbb{Z}_{\ell_1^{v_{\ell_1}(N)}} \oplus \mathbb{Z}_{\ell_2^{v_{\ell_2}(N)}} \oplus \cdots \oplus \mathbb{Z}_{\ell_k^{v_{\ell_k}(N)}},$$

($\ell \ne p$) with $r_\ell + s_\ell = v_\ell(N)$ and $\min(r_\ell, s_\ell) \le v_\ell(q-1)$.

Some specific details on the actual decomposition of the groups of points as the product of two cyclic groups are as follows.

***Theorem 46.*** ([CJ02]). Let $E$ be a CM elliptic curve, and let $a^2 - \Delta b^2 = 4q$. Then $\#E(\mathbf{F}_q) \cong \mathbb{Z}_d \times \mathbb{Z}_{de}$, where $d = \gcd(b, (a + b\delta - 2)/2)$, $e = (p + 1 - a)/d^2$, and

$$\delta = \begin{cases} 0 & \Delta \equiv 0 \bmod 4, \\ 1 & \Delta \equiv 1 \bmod 4. \end{cases} \quad (35)$$

From the expression $N = q + 1 - t = d^2 e \le q + 1 + 2q^{1/2}$, it is clear that $d \le 2q^{1/2}$, and $e \ge q^{1/2}$. Equality (maximal or minimal curves), seems to occur only with some supersingular elliptic curves. Moreover, from Cassel's Theorem, it follows that $E(\mathbf{F}_q)$ is cyclic if $N = \#E(\mathbf{F}_q)$ is squarefree.

***Proposition 47.*** The group $E(\mathbf{F}_q) \cong \mathbb{Z}_d \times \mathbb{Z}_{de}$ if and only if there are two independent points of order $d$ in $E(\mathbf{F}_q)$. Two independent points occur if and only if $d$ divides both $q - 1$, and $q + 1 - t$. Specifically, if $q + 1 - t = d^2 e$, and $q - 1 = df$.

The *exponent* $e$ of an elliptic curve $E$ is the size of the largest cyclic subgroup $\mathbb{Z}_e$ contained in $E(\mathbf{F}_q)$.

***Theorem 48.*** ([DK00]). Let $E$ be an elliptic curve defined over $\mathbb{Q}$. If $E$ does not have CM, assume the GRH. Let $f(x)$ be any positive function on $[2, \infty)$ that tends to infinity with $x$. Then the exponent $e$ of $E$ satisfies $e > p^{3/4}/\log p$ for almost all $p$.

*Miller Algorithm*

The Miller algorithm determines the group structure of $E(\mathbf{F}_q)$ as either a cyclic group $E(\mathbf{F}_q) \cong \mathbb{Z}_N$ or as the product of two cyclic groups $E(\mathbf{F}_q) \cong \mathbb{Z}_d \times \mathbb{Z}_{de}$. It requires the complete factorization of $\gcd(q - 1, \#E(\mathbf{F}_q))$, and runs in random polynomial time.

*Example* 49.  This illustrates the variety of groups $E(\mathbf{F}_{13})$ of elliptic curves over the finite field $\mathbf{F}_{13}$, and how to determine the group structure.

(i) The group of points of the elliptic curve $E : y^2 = x^3 + 2$ over $\mathbf{F}_{13}$ has order $19 = \#E(\mathbf{F}_{13})$, see Example 6. Since $p - 1 = 2^2 \cdot 3$ and $p + 1 - t = d^2 e = 19$, the parameters are $d = 1$, $e = 19$. Thus, the group is cyclic $E(\mathbf{F}_q) \cong \mathbb{Z}_{19}$. In fact, any points generates the group, viz, $E(\mathbf{F}_q) = <(1, 4)>$.

(ii) The group of points of the twisted curve $E' : y^2 = x^3 + 3$ is $E'(\mathbf{F}_{13}) = \{ (\infty, \infty), (0, 1), (0, 9), (0, 4), (1, 2), (1, 11), (3, 11), (9, 2), (9, 11) \}$ has cardinality $9 = \#E'(\mathbf{F}_{13})$. Since $p - 1 = 2^2 \cdot 3$ and $p + 1 - t = d^2 e = 9$, the parameters are $d = 3$, $e = 1$. Here, $d = 3$ divides both $p - 1$ and $p + 1 - t$. Thus, the group is not cyclic but the product of two cyclic $E'(\mathbf{F}_q) \cong \mathbb{Z}_3 \times \mathbb{Z}_3$. In fact, it is generated by two points of order 3 each, viz, $E'(\mathbf{F}_q) = <(1, 2), (0, 9)>$.

(iii) The group of the nonCM curve $E'' : y^2 = x^3 + 6x + 11$ is $E''(\mathbf{F}_{13}) = \{ (\infty, \infty), (3, \pm 2), (5, \pm 6), (8, \pm 5), (9, \pm 1) (12, \pm 2) \}$ has order $15 = \#E''(\mathbf{F}_{13})$. Since $p - 1 = 2^2 \cdot 3$ and $p + 1 - t = d^2 e = 15$, the parameters are $d = 1$, $e = 15$. Thus, the group is cyclic $E''(\mathbf{F}_q) \cong \mathbb{Z}_{15} \cong \mathbb{Z}_3 \times \mathbb{Z}_5$. In fact, it is generated by a single point, viz, $E''(\mathbf{F}_q) = <(5, 7)>$.

## 21. *n*-Torsion Group $E(\mathbf{F}_q)[n]$

Let $n \in \mathbb{N}$, the *n*-torsion group $E(\mathbf{F}_q)[n]$ ( or $E[n]$ for short) is the set of points $P \in E(\mathbf{F}_q)$ such that $nP = O$. The *n*-torsion is the kernel of the multiplication map. More precisely, $E[n] = \ker([n]) = \{ P \in E(\mathbf{F}_q) : [n]P = O \}$. It comes in two varieties.

***Theorem* 50.**  (Tate 74).  The group of points $E(\overline{\mathbf{F}}_q)$ is a torsion group such that

1. $E(\overline{\mathbf{F}}_q)[n] \cong \mathbb{Z}_n \times \mathbb{Z}_n$,     if $p$ does not divide $n$, and $q = p^e$,
2. $E(\overline{\mathbf{F}}_q)[n] \cong \mathbb{Z}_n$,     if $p$ divides $n$, and the elliptic curve is not supersingular,
3. $E(\overline{\mathbf{F}}_q)[n] \cong \{ O \}$     if $p$ divides $n$, and the elliptic curve is supersingular.

Recall that $E(\mathbf{F}_q)$ is the kernel of the Frobenius $\sigma - 1$ action on the group $E(\overline{\mathbf{F}}_q)$. Note: The kernel ker($f$) of a function $f$ over a group is the subset of elements that are mapped to the identity, for example, the kernel of $f(x) = x^2$ over the multiplicative group of real numbers $\mathbb{R}^*$ is ker($f$) = $\{ x \in \mathbb{R}^* : x^2 = 1 \} = \{ -1, 1 \}$.

**2-Sylow Subgroup of $E(\mathbf{F}_q)$**
The structure of the 2-Sylow subgroup $S_2(\mathbf{F}_q) = E(\mathbf{F}_q)[2^n]$ of $E(\mathbf{F}_q)$ is determined by the number of roots of $f(x) = x^3 + ax^2 + bx + c$ in $\mathbf{F}_q$.

(i) If the cubic polynomial $f(x)$ has a single root, then subgroup $E(\mathbf{F}_q)[2] \cong \mathbb{Z}_2$, which implies that $E(\mathbf{F}_q)[2^n]$ is cyclic of rank 1.

(ii) If the cubic polynomial $f(x)$ has three roots, then the subgroup $E(\mathbf{F}_q)[2] \cong \mathbb{Z}_2 \times \mathbb{Z}_2$, which implies that $E(\mathbf{F}_q)[2^n]$ is not cyclic and of rank 2.

An algorithm for computing $S_2(\mathbf{F}_q)$ based on successive point halving is given in [MT05]. Two parametized infinite families of curves whose 2-Sylow subgroups are either cyclic or noncyclic are also described.

## 22. Inverse Points Counting

The inverse point counting problem starts with a prescribed cardinality $N > 1$ and produces an algebraic curve such that $N = \#E(\mathbf{F}_q)$. The inverse counting problem of some orders $N = \#E(\mathbf{F}_q)$ has an easy solution:

(i) $N = q$,    anomalous elliptic curves,
(ii) $N = q + 1$,   supersingular elliptic curves,

and some other cases. But the general case for arbitrary orders $N$ is much more difficult, and it is handled using several techniques.

A random elliptic curve whose group of points $E(\mathbf{F}_q)$ has order $N \in [q + 1 - 2q^{1/2}, q + 1 + 2q^{1/2}]$ can be determined in probabilistic exponential time $O(N^{1/2+\varepsilon})$, $\varepsilon > 0$, see [BR04]. A sketch of the procedure is described below. This procedure assumes that any pair of integers $N$, and $q$ is admissible in the construction of a curve.

*Algorithm* I.
Input: $N$, and $q$.
Output: $E : y^2 = x^3 + ax + b$ of order $N$ over $\mathbf{F}_q$.
1. Select a random $0 \ne a \in \mathbf{F}_q$, and put $E : y^2 = x^3 + ax - a$.
2. Compute the scalar multiple $NP = ?$ of the point $P = (1, 1)$.
3. If $NP = O$, then compute the order $\#E(\mathbf{F}_q)$. Otherwise repeat step 1.
4. If $\#E(\mathbf{F}_q) = N$, then return $E$, else return the twisted elliptic curve $E' : y^2 = x^3 + as^2 x + bs$, where s is a quadratic nonresidue.

## 23. Complex Multiplication Method

The complex multiplication method is a mean of constructing curves with groups of points of given orders. Given a prime power and $N \in [q + 1 - 2q^{1/2}, q + 1 + 2q^{1/2}]$, the algorithm finds the algebraic curve. Basically, an integral solution of the equation $4qg^2 = t^2 + du^2$, $q$ fixed, and $d$ a discriminant, e.g., $d \equiv 3 \mod 4$, is converted into an algebraic curve of genus $g > 0$, whose group of points has the cardinality $N = \#E(\mathbf{F}_q) = q + 1 - t$.

The characteristic polynomial of the Frobenious endomorphism is $L(x) = x^2 - tx + q$, and its value at $x = 1$ gives $L(1) = \#E(\mathbf{F}_q) = q + 1 - t$.

**Lemma 51.**  If a nonsingular algebraic curve of genus $g > 0$ has $N = q + 1 - t$ points, then $4qg^2 = t^2 + du^2$, for some $d \ge 1$.

Proof: By Hasse's Theorem $|t| \le 2g\sqrt{q}$. Thus $4qg^2 - t^2 = v = du^2 \ge 0$.    ∎

**Definition 52.**  A squarefree integer $d$ is a *fundamental discriminant* if $d \equiv 1 \mod 4$ or $d \equiv 8, 12 \mod 16$.

If $d$ is a fundamental discriminant, than the set $\text{End}(E)$ is the full ring of integers in $\mathbb{Q}(\sqrt{-d})$, otherwise, it is an order strictly contained in $\mathbb{Q}(\sqrt{-d})$.

If the prime power $q$ splits completely in the quadratic field of $\mathbb{Q}(\sqrt{d})$ of discriminant $d < 0$, then there is an algebraic integer $\pi \in \mathbb{Q}(\sqrt{d})$ such that $4q = t^2 - du^2 = \pi\bar{\pi}$, where $\pi = t + u\sqrt{d}$.

Moreover, if $d \equiv 1 \bmod 4$, or $8, 12 \bmod 16$ is a fundamental discriminant and $u = 1$, then the ring $\text{End}(E) = \mathbb{Q}(\sqrt{d})$. The integer $u$ is the index $[\mathbb{Q}(\sqrt{d}) : \text{End}(E)] = u \geq 1$ of $\text{End}(E)$ in $\mathbb{Q}(\sqrt{d})$. The curve is not supersingular if $t \neq 0$ and its $\text{End}(E) \subseteq \mathbb{Q}(\sqrt{d})$.

**Theorem 53.** ([LA94]). Let $E$ be an elliptic curve such that $E(\mathbf{F}_q) = \mathbb{Z}_m \times \mathbb{Z}_n$. Then one of the followings holds
1. The ring $\text{End}(E)$ is an order in $\mathbb{Q}(\sqrt{-1})$ and with $m \mid n$, and $p = n^2 + 1$.
2. The ring $\text{End}(E)$ is an order in $\mathbb{Q}(\sqrt{-3})$ and with $m \mid n$, and $p = n^2 \pm n + 1$.

**Example 54.** (i) The of elliptic curve $E : y^2 = x^3 + b$ has complex multiplication by $\mathbb{Q}(\sqrt{-1})$, so $\text{End}(E)$ contains $\mathbb{Q}(\sqrt{-1})$. This implies that $E(\mathbf{F}_p) \cong \mathbb{Z}_n \times \mathbb{Z}_n$ for any prime $p = n^2 + 1$ with $\pi = \pm 1 \pm n\sqrt{-1}$. Observe that the trace $t = Tr(\pi) = (\pi + \bar{\pi})/2 = \pm 1$ remains constant as the prime $p$ varies.

(ii) The of elliptic curve $E : y^2 = x^3 + ax$ has complex multiplication by $\mathbb{Q}(\sqrt{-3})$, so $\text{End}(E)$ contains $\mathbb{Q}(\sqrt{-3})$. This implies that $E(\mathbf{F}_p) \cong \mathbb{Z}_n \times \mathbb{Z}_n$ for any prime $p = n^2 \pm n + 1$ with $\pi = (\pm n - 2 \pm n\sqrt{-3})/2$. Observe that the trace $t = Tr(\pi) = (\pi + \bar{\pi})/2 = \pm n - 2$ does not remain constant as the prime $p$ varies.

*Hilbert Class Polynomial*
Hilbert class polynomial is defined by

$$H_d(x) = \prod_{(a,b,c)} (x - j((-b + \sqrt{d})/2a) \in \mathbb{Z}[x], \tag{36}$$

where the index $(a, b, c)$ ranges over the $h = h(d)$ reduced quadratic forms $ax^2 + bxy + cy^2$ of discriminant $d < 0$. The integer $h(d)$ is the class number of the quadratic field $\mathbb{Q}(\sqrt{d})$.

A root $j$ of $H_d(x)$ is the $j$-invariant of an isomorphism class of elliptic curves over $\mathbf{F}_q$ with $\text{End}(E) \subseteq \mathbb{Q}(\sqrt{d})$. The basic pattern of the class polynomials and the representative elliptic curves in characteristic $p \neq 2, 3$ are as follows:

The $j$-invariant of an elliptic curve over the rational numbers with complex multiplication is an integer. There are nine cases of integral $j$-invariants associated with the quadratic fields of class number one, $h(d) = 1$. For these determined nine cases the polynomial $H_d(x)$ are linear:

| D | $\tau$ | $j(\tau)$ | $H_d(x)$ |
| --- | --- | --- | --- |
| $-3$ | $(1 + \sqrt{-3})/2$ | $0$ | $X$ |
| $-4$ | $\sqrt{-1}$ | $12^3$ | $x - 12^3$ |
| $-7$ | $(3 + \sqrt{-7})/2$ | $-15^3$ | $x + 15^3$ |
| $-8$ | $\sqrt{-2}$ | $20^3$ | $x - 20^3$ |
| $-11$ | $(3 + \sqrt{-11})/2$ | $-32^3$ | $x + 32^3$ |
| $-19$ | $(3 + \sqrt{-19})/2$ | $-96^3$ | $x + 96^3$ |
| $-43$ | $(3 + \sqrt{-43})/2$ | $-960^3$ | $x + 960^3$ |
| $-67$ | $(3 + \sqrt{-67})/2$ | $-5280^3$ | $x + 5280^3$ |
| $-163$ | $(3 + \sqrt{-163})/2$ | $-640320^3$ | $x + 640320^3$ |

(1) $H_1(x) = x - 12^3$ for the isomorphism class of elliptic curves $E : y^2 = x^3 + ax$, $a \neq 0$, over $\mathbf{F}_q$ of $j$-invariant $j(E) = 12^3$. There are four inequivalent twisted elliptic curves $E' : sy^2 = x^3 + ax$ by $s \in \mu_4 = \{\pm 1, \pm i\}$.

(2) $H_3(x) = x$ for the isomorphism class of elliptic curves $E : y^2 = x^3 + b$, $b \neq 0$, over $\mathbf{F}_q$ of $j$-invariant $j(E) = 0$. There are six inequivalent twisted elliptic curves $E' : sy^2 = x^3 + b$ by $s \in \mu_6 = \{1, \omega, \omega^2, \ldots, \omega^5\}$.

(3) For $H_d(x)$, $d = -7, -8, -11, -19, -43, -67$, and $-163$, for the isomorphism classes of seven elliptic curves $E : y^2 = x^3 + ax + b$, $b \neq 0$, over $\mathbf{F}_q$ of $j$-invariant $j(E) \neq 0, 12^3$. There are two inequivalent twisted elliptic curves $E' : sy^2 = x^3 + b$ by $s \in \mu_2 = \{\pm 1\}$.

For $h = h(d) > 1$, the polynomial $H_d(x) = a_h x^h + \cdots + a_1 x + a_0$ is associated to the isomorphism class of elliptic curves $E : y^2 = x^3 + 3ax + 2a$, $a = j/(1728 - j)$ over $\mathbf{F}_q$ of $j$-invariant $j(E) = j \neq 0, 1728$. This curve has two twisted curves by $s \in \mu_2 = \{\pm 1\}$.

The difficult part of this procedure is the calculation of the polynomial $H_d(x)$. The complete time complexity analysis of the algorithm for computing the Hilbert class polynomial, in about $O(d(\log q)^3)$ steps, is given in [AG04], see also [LA94].

**Anomalous Curves**

An elliptic curve $E : y^2 = x^3 + ax + b$ over $\mathbf{F}_q$ is called *anomalous* if $\#E(\mathbf{F}_q) = q$ and its twist $E' : y^2 = x^3 + ac^2 x + bc^3$ has $\#E'(\mathbf{F}_q) = q + 2$ points.

Anomalous elliptic curves can be viewed as having embedding degrees $k = 0$. For example, the projection map $(x, y) \to y$ is a one-to-one, which maps $E(\mathbf{F}_q)$ to $\mathbf{F}_q$.

Let $d \equiv 1 \bmod 4$, $d < 0$, and choose a prime $p$ such that $4p = 1 - du^2$. Then the CM method produces a polynomial whose roots are the $j$-invariants of elliptic curves with the endomorphisms ring $\text{End}(E)$ contained in $\mathbb{Q}(\sqrt{d})$. Here, the trace is $t = \alpha + \bar{\alpha} = (\pm 1 - u\sqrt{d})/2 + (\pm 1 + u\sqrt{d})/2 = \pm 1$. Thus, the curve $E$ or its twist $E'$ has $p$ $\mathbf{F}_p$-rational points. Thus, its twisted curve is an anomalous curve.

## 24. Linear Pairing

Let $n \geq 1$ be a large integer, say prime, and consider the minimal extension $\mathbf{F}_{q^k}$ of the finite field $\mathbf{F}_q$, which has the smallest embedding degree $k$ such that $n \mid q^k - 1$, but $n \nmid q^d - 1$ for all $d < k$. Let $G$ and $H$ be subgroups of rational points $E(\mathbf{F}_{q^k})$, where the order $\#G = n$. A linear pairing is a function $< , > : G \times H \to \mathbf{F}_{q^k}$. There are several varieties of linear pairings. The most popular are the Weil pairing and Tate pairing.

**Weil Pairing.** Let $E(\mathbf{F}_{q^k})[n] = \{P \in E(\mathbf{F}_{q^k}) : [n]P = O\}$ be the subset of $n$-torsion points. The Weil pairing maps a point in cross product $w : E(\mathbf{F}_{q^k})[n] \times E(\mathbf{F}_{q^k})[n] \to \mu_n \subset \mathbf{F}_{q^k}$ to a root of unity in the finite field extension $\mathbf{F}_{q^k}$.

***Properties of Weil Pairing.*** Let $P, Q \in E(\mathbf{F}_{q^k})$ be a pair of points. Then, the following properties hold:

1. $e(P, O) = e(O, P) = 1$,  Identity.
2. $e(P, Q) = e(Q, P)^{-1}$,  Antisymmetry, or noncommutative.
3) $e(P, P) = 1$,  Trivial.
4. $e(P, Q) \neq 1$ for some $P, Q \in E(\mathbf{F}_{q^k})$,  Nonanomalous.
5. $e(P + Q, R) = e(P, R)e(Q, R)$,  Bilinearity.

6. $e(-P, Q) = e(P, -Q) = e(P, Q)^{-1}$,                  Inverse.

**Tate Pairing.** Let $E(\mathbf{F}_{q^k})[n] = \{ P \in E(\mathbf{F}_{q^k}) : [n]P = O \}$ be the subset of $n$-torsion points, and let $nE(\mathbf{F}_{q^k}) = \{ nP : P \in E(\mathbf{F}_{q^k}) \}$ be a coset of points. The Tate pairing maps a point in the cross product $w : E(\mathbf{F}_{q^k})[n] \times nE(\mathbf{F}_{q^k}) \to \mathbf{F}_{q^k}/(\mathbf{F}_{q^k})^n \subset \mathbf{F}_{q^k}$ to the quotient. Since $\mathbf{F}_{q^k}/(\mathbf{F}_{q^k})^n = \{ x = yz^n : x, y, z \in \mathbf{F}_{q^k} \}$, the image of a point $P \in E(\mathbf{F}_{q^k})$ is unique up to an $n$th root of unity in $\mu_n \subset \mathbf{F}_{q^k}$.

***Properties of Tate Pairing.*** Let $P, Q \in E(\mathbf{F}_{q^k})$ be a pair of points. Then, the following properties hold:

1. $w(P, O) = w(O, P) = 1$,                                            Identity.
2. $w(P, Q) = w(Q, P)$,                                     Symmetry, or noncommutative.
3) $w(P, P) \neq 1$,                                                 Nontrivial.
4. $w(P, Q) \neq 1$ for some $P, Q \in E(\mathbf{F}_{q^k})$,            Nonanomalous.
5. $w(P + Q, R) = w(P, R)w(Q, R)$,                    Bilinearity.
6. $w(-P, Q) = w(P, -Q) = w(P, Q)^{-1}$,             Inverse.

The verification of these properties, for both Weil and Tate pairings, are easy consequences of the cyclic structure of the group of $n$th root of unity $\mu_n$, see [GS02], [BL05] for implementations ideas, and [FR94], [SN92] et cetera, for the theory of pairings.

An important application of the Weil and Tate pairings are the transfer of the discrete logarithm on elliptic curves to the discrete logarithm on finite fields. This is accomplished by means of the Weil pairing or the Tate pairing

$$e(aP, bQ) = e(P, Q)^{ab} \quad \text{or} \quad w(aP, bQ) = w(P, Q)^{ab}$$

of the points $aP, bQ \in E(\mathbf{F}_{q^k})$, where $a, b \in \mathbb{Z}$ are integers. In the case of small embedding degree $k$, it effectively reduces the complexity from fully exponential time to subexponential time, see [MV93], [FR94] and similar references.

**Small Embedding Degrees.** A random elliptic curve has negligible probability of having small embedding degree $k = O((\log q)^c)$, $c > 0$ constant, see [BK98]. However, there are methods for constructing elliptic curves of small embedding degrees $k$.

***Lemma 55.*** Let $r$ be a prime divisor of $\#E(\mathbf{F}_q)$ and $k > 0$ be its embedding degree. Then
(i) $t - 1$ is a $k$th root of unity in $\mathbf{F}_r$.          (ii) $t \equiv \omega + 1 \bmod r$, $\omega$ is a $k$th root of unity in $\mathbf{F}_r$.

These simple observations are used to construct a curve of the given orders, see [BR03].

***Definition 56.*** Let $r$ be a divisor of $\#E(\mathbf{F}_q)$. The embedding degree of the group of points $E(\mathbf{F}_q)$ is the smallest integer $k > 0$, such that
(i) $N = q + 1 - t \equiv 0 \bmod r$,                       (ii) $q^k - 1 \equiv 0 \bmod r$.

***Theorem 57.*** ([KZ98]). Let $E$ be an elliptic curve defined over $\mathbf{F}_q$, and suppose that $r$ is a prime that divides $N = E(\mathbf{F}_q)$ but does not divide $q - 1$. Then $E(\mathbf{F}_{q^k})$ contains $r^2$ points of order $r$ if and only if $r \mid q^k - 1$.

***Proposition 58.*** ([MV93]) The class of supersingular elliptic curves has embedding degrees of $k = 1, 2, 3, 4,$ and 6.

For each $k < 7$, there is a simple proof based on divisibility conditions. For example, to prove the case $k = 6$, observe that

$$\#E(\mathbf{F}_q) = q+1\pm\sqrt{3q} \Rightarrow (q+1)(q+1-\sqrt{3q})(q+1+\sqrt{3q}) = q^3+1 \Rightarrow k=6$$

since $(q^3 + 1)(q^3 - 1) = q^6 - 1$ is the smallest integers such that $\#E(\mathbf{F}_q) = q^3 + 1$ divides $q^k - 1$, see [BL05, p.199] for more details.

A list of low embedding degree $k$ and the group structure of the corresponding elliptic curves is included here.

| K | Q | $\#E(\mathbf{F}_q)$ | Structure of $E(\mathbf{F}_{q^k}) = \mathbb{Z}_m \times \mathbb{Z}_n$ |
|---|---|---|---|
| 1 | $p^{2m}$ | $q+1\pm 2\sqrt{q}$ | $q^{1/2}\pm 1, q^{1/2}\pm 1$ |
| 2 | $p^{2m+1}, p^{2m}$ with $p = 4a + 3$ | $q+1$ | $q+1, q+1$ |
| 3 | $p^{2m}$ with $p = 3a + 2$ | $q+1\pm\sqrt{q}$ | $q^{3/2}\pm 1, q^{3/2}\pm 1$ |
| 4 | $2^{2m+1}$ | $q+1\pm\sqrt{2q}$ | $q^2\pm 1, q^2\pm 1$ |
| 6 | $3^{2m+1}$ | $q+1\pm\sqrt{3q}$ | $q^3+1, q^3+1$ |

**Proposition 59.** Let $E_0 : y^2 = x^3 + 2x + 1$ and $E_1 : y^2 = x^3 + 2x - 1$ be curves over $\mathbf{F}_3$. Then

$$\#E_0(\mathbf{F}_{3^n}) = \begin{cases} 3^n + 1 + \sqrt{3^{n+1}} & n \equiv \pm 1 \bmod 12, \\ 3^n + 1 - \sqrt{3^{n+1}} & n \equiv \pm 5 \bmod 12, \end{cases} \tag{37}$$

and

$$\#E_1(\mathbf{F}_{3^n}) = \begin{cases} 3^n + 1 - \sqrt{3^{n+1}} & n \equiv \pm 1 \bmod 12, \\ 3^n + 1 + \sqrt{3^{n+1}} & n \equiv \pm 5 \bmod 12. \end{cases}$$

These supersingular curves have embedding degree $k = 6$, so the discrete logarithms on these curves can be transferred to the finite field, that is, $E_i(\mathbf{F}_{3^n}) \to \mathbf{F}_{3^{6n}}$ by means of pairing, see [MV93].

The discrete logarithm in groups of embedding degree $k < (\log q)^2$ has subexponential time complexity, this follows from the time complexity $\exp(c(\log q^k)^{1/3}(\log\log q^k)^{2/3})$, $c > 0$, of the discrete logarithm in finite fields.

*The Distortion Map.* The distortion map $(x, y) \to \psi(x, y) = (\rho x, \pm y)$ exchanges the point $P = (x, y) \in E(\mathbf{F}_q)$ for a different point $\psi(P) = (\rho x, \pm y)$.

The distortion map is used in the construction of groups isomorphisms. Specifically, the isomorphism of a subgroup $G = <P>$ of $E(\mathbf{F}_q)$ of order $n = \#G > 2$ to a subgroup $\mu_n = \{x \in \mathbf{F}_{q^k} : x^n - 1 = 0\}$ of root of unity in $\mathbf{F}_{q^k}$. There are other applications of this map.

| k | Char $p$ | Aut($E$) | Curve | Distortion Map | $\#E(\mathbf{F}_q)$ |
|---|---|---|---|---|---|
|   | $p = 2$ | $j = 0$ | $??y^2 + y = x^3$ | | |
|   | $p = 2$ | $j \neq 0$ | $??y^2 + xy = x^3 + x^2 + j^{-1}$ | | |
|   | $p = 3$ | | $y^2 = x^3 + x + c$ | | $3^n + 1 \pm 3^{(n+1)/2}$ |
| 2 | $p \equiv 2 \bmod 3$ | $\{1, \rho, \rho^2\}$ | $y^2 = x^3 + ax, a \neq 0$ | $\psi(x, y) = (\rho x, \pm y)$ | $p + 1$ |
| 2 | $p \equiv 3 \bmod 4$ | $\{\pm 1, \pm i\}$ | $y^2 = x^3 + b, b \neq 0$ | $\psi(x, y) = (-x, iy)$ | $p + 1$ |

## 25. Cyclicity of $E(\mathbf{F}_p)$

A squarefree order $\#E(\mathbf{F}_q)$ is a sufficient condition to have a cyclic group $E(\mathbf{F}_q) \cong \mathbb{Z}_N$. , this follows from Cassel's Theorem. Moreover, since the number of squarefree integers $N < q + 1 + 2q^{1/2}$ is asymptotic to $(6/\pi^2)q$, the set of elliptic curves that have cyclic groups is large.

***Theorem* 60.** ([CJ03]). Let $E$ be an elliptic curve over $\mathbb{Q}$ of conductor $N_0$ and with complex multiplication by the full ring of integers of $\mathbb{Q}(\sqrt{-d}\,)$. Then the following hold.
(i) As $x \to \infty$

$$f(x, \mathbb{Q}) = \#\{\, p \leq x : p \nmid N_0 \text{ and } E(\mathbf{F}_p) \text{ is cyclic } \} = c_E \mathrm{li}(x) + O\left(\frac{x}{\log x \log \log \log(x)}\right), \tag{38}$$

where $\mathrm{li}(x)$ is the logarithmic integral, and $c_E > 0$ is a constant and the implied constant depends on $N$.

(ii) The smallest prime $p \nmid N_0$ for which $E(\mathbf{F}_p)$ is cyclic has size $O(e^{N^2})$.

Under the assumption of the generalized Riemann hypothesis for the zeta function of an elliptic curve over $\mathbb{Q}$ of conductor $N_0$, and such that $\mathbb{Q}(E[2]) \neq \mathbb{Q}$, a sharper estimate is known, details are given in [CM04]:
(iii) If $E$ has CM by the full ring of integers, then the smallest prime $p \nmid N_0$ for which $E(\mathbf{F}_p)$ is cyclic has size $p = O(\log(N)^{2+\varepsilon})$, $\varepsilon > 0$.
(iv) If $E$ is a nonCM curve, then the smallest prime $p \nmid N_0$ for which $E(\mathbf{F}_p)$ is cyclic has size $p = O(\log(N_0)^{4+\varepsilon})$, $\varepsilon > 0$.

## 26. Distribution of Prime and Nearly Prime Cardinalities

Cryptographic consideration calls for groups of prime or nearly prime orders. The primality of $\#E(\mathbf{F}_q)$ demands that the polynomial $f(x)$ in $E : y^2 = f(x)$ has no roots in the finite field $\mathbf{F}_q$ of definition. This condition is often denoted by $\mathbb{Q}(E[2]) \neq \mathbb{Q}$.

Let $B > 0$. An integer $N > 0$ is called $B$-smooth if all of its prime factors are less than $B$. A question of much interest in elliptic curve cryptography and elliptic curve method of integer factorization is the smoothness of the cardinality $N = \#E(\mathbf{F}_q)$ of the groups of points of a random curve. In elliptic curve cryptography there is an interest in having nonsmooth cardinalities $N = \#E(\mathbf{F}_q)$ to avoid several well-known effective discrete logarithms algorithms on groups of smooth cardinalities. For example, the order of a subgroup $G$ of $E(\mathbf{F}_q)$ should be divisible by a large prime.

The number of prime factors $\omega(N)$ in $N = \#E(\mathbf{F}_q)$ falls in the Ramanujan interval

$$\log\log(N) - \log\log(N)^{1/2+\varepsilon} \leq \omega(N) \leq \log\log(N) + \log\log(N)^{1/2+\varepsilon}$$

for almost all prime powers $q$, and $\varepsilon > 0$. The exceptions consist of a subset of prime powers of zero density in the set of prime powers. Therefore, most group of points $E(\mathbf{F}_q)$ have composite orders, but have very few prime factors. Also observe that due to the divisibility property $N_d \mid N_n$ if $d \mid n$, the groups $E(\mathbf{F}_{p^n})$ = prime must have $n$ prime.

The number of primes $p$ for which the cardinality $N = \#E(\mathbf{F}_p)$ is prime is not known, but it is believed to have the asymptotic behavior

$$\#\{ p \leq x : \#E(\mathbf{F}_p) \text{ is prime} \} \sim \frac{cx}{\log(x)^2}. \tag{39}$$

The best proven results state the following. Here, the number theoretical function $\omega(N) = \#\{ p \mid N : p \text{ is prime} \}$, and $\Omega(N) = \#\{ p^i \mid N : p^i \text{ is a prime power} \}$.

**Theorem 61.** ([ST05]). Assuming the extended Riemann hypothesis, and that $E$ has finitely many isogenous curves. Let $N = \#E(\mathbf{F}_p)$. Then

(i) If $E$ has complex multiplication then $\#\{ p \leq x : \Omega(N) \leq 3 \} \geq \dfrac{cx}{\log(x)^2}$.

(ii) If $E$ does not have complex multiplication then $\#\{ p \leq x : \omega(N) \leq 5 \} \geq \dfrac{cx}{\log(x)^2}$.

A proof conditioned on a weak generalized Riemann hypothesis for nonCM elliptic curves is carried out in [CJ05].

**Theorem 62.** ([CJ05]). If $E$ is a nonCM elliptic curve of conductor $N_0$ over the rational, then

$$\#\{ p \leq x : p \nmid N \text{ and } \#E(\mathbf{F}_q) = \text{prime} \} \leq \frac{cx}{\log(x) \log \log \log(x)}. \tag{40}$$

This should be compared to the conjectured size of this set which is $O(x/\log(x)^2)$.

## 27. Primitive Points on Curves

**Theorem 63.** ([SF95]). Let $p > 229$ be a prime and let $E$ be an elliptic curve over $\mathbf{F}_p$. Then either $E$ or its quadratic twist $E'$ admits an $\mathbf{F}_p$-rational point $P$ with the property that the only integer $n \in [q+1-2\sqrt{q}, q+1+2\sqrt{q}]$ for which $nP = O$ is the order of the group of points.

It readily follows that either the elliptic curve $E$ or its quadratic twist $E'$ has a cyclic group of points. A primitive point that generates the entire group of points can be determined by a brute force search.
The worst case uses up to $O(q^{3/2})$ elliptic curves operations. A better algorithm that runs in random exponential time $O(q^{3/2})$ is sketched here.

*Primitive Points Algorithm* I
Input: A nonsingular curve $E : y^2 = f(x)$.
Output: A primitive point $P$.
1. Compute $n = \#E(\mathbf{F}_q)$.
2. Select a small set of random points $P_0, P_1, \ldots, P_d \in E(\mathbf{F}_q)$, some $d \leq \log q$.
3. Compute the orders $n_0, n_1, \ldots, n_k \leq q + 1 + 2q^{1/2}$ of the points such that $n_0 P_0 = O, n_1 P_1 = O, \ldots, n_k P_k = O$, and $\gcd(n_i, n_j) = 1$.
4. If $n = \text{lcm} \{ n_0, n_1, \ldots, n_k : \gcd(n_i, n_j) = 1 \} > q + 1 - 2q^{1/2}$, then return the primitive point $P_0 + P_1 + \cdots + P_k$ of order $n = \#E(\mathbf{F}_q)$.
5. Repeat steps 2 to 4 for the twisted curve $E' : sy^2 = f(x)$, $0 \neq s$ a quadratic nonresidue.

The core of this technique is essentially the Gaussian algorithm for computing primitive roots in finite fields. This is coupled to the effective range of the order specified by Hasse's theorem, which gives the time complexity of the algorithm.

## 28. Distribution of the Angles of $a_p$

Let $p$ be a prime, let $E: y^2 = x^3 + ax + b$ be a fixed elliptic curve over $\mathbf{F}_p$ and let $a_p = 2\sqrt{p}\cos\theta$ be its trace. It is the claim of a long standing problem that the normalized angle $\theta$ is uniformly distributed in the interval $[-\pi, \pi]$ with respect to the probability measure $\mu = 2\pi^{-1}\sin^2(x)dx$ as prime $p$ varies.

**Sato-Tate Conjecture 64.** Let $\pi(x) = \#\{$ prime $p \leq x \}$, and let $a_p = 2\sqrt{p}\cos\theta$. Then

$$\lim_{x\to\infty} \frac{1}{\pi(x)} \#\{ p \leq x : a \leq \theta \leq b \} = \frac{2}{\pi}\int_a^b \sin^2(t)dt. \tag{41}$$

This conjecture has been settled for the CM case, but it is not settled for the nonCM case.

**Theorem 65.** Let $E$ be CM elliptic curve. Then the angles are uniformly distributed on $[0, \pi]$ with respect to the probability measure $d\mu = (\pi\delta_{\pi/2} + d\theta)/2$. Specifically,

$$\lim_{x\to\infty} \frac{1}{\pi(x)} \#\{ p \leq x : a \leq \theta \leq b \} = \frac{1}{2\pi}\left(\pi\delta_{\pi/2}I_{[a,b]} + b - a\right), \tag{42}$$

where $I_{[a,b]}(x)$ is the characteristic function of the interval $[a, b]$.

This is derived from the work of Hecke on the uniform distribution of the pair of angles $(-\theta, \theta)$, see [MY82], [ML05] for the proof, and discussions of the Sato-Tate conjecture, which is a proven Theorem over function fields.

The distributions of the angles of a handful of other exponential sums and complex-valued quantities have been proven or are believed to obey similar uniform distribution laws. For example, Gaussian sums, Kloosterman sums, Gaussian primes, etc.

For curve $C$ of fixed genus $g > 0$, and $\#C(\mathbf{F}_q) = q + 1 - a_q$ such that $|a_q| \leq 2g\sqrt{q}$, the random variable $a_q/(2\sqrt{q})$ is assumed to be uniformly distributed in $[-g, g]$ with respect to the probability measure $\mu_q = 2\pi^{-1}\sin^2(x)dx$.
The extension of the Sato-Tate measure to random curves of random genus $g > 0$ is considered in [LN02]. In this setting the angle has a (asymptotically) normal distribution. More precisely, for random curve $C : y^2 = f(x)$ of random genus $g > 0$, and $\#C(\mathbf{F}_q) = q + 1 - a_q$ such that $|a_q| \leq 2g\sqrt{q}$, the random variable $a_q/(2\sqrt{q})$ is normally distributed in $[-\infty, \infty]$ with respect to the measure $\mu = (2\pi)^{-1/2} e^{-x^2/2} dx$.

**Theorem 66.** Fix a prime number $p$ and fix an elliptic curve $E_n : y^2 = x^3 + ax + b$ over $\mathbf{F}_p$. Let $N_n = p^n + 1 - a_n$ If the trace $a_1 = 2\sqrt{p}\cos\theta_1$, where $\theta_1 \neq r\pi$ is irrational, $r \in \mathbb{Q}$, then the angle $\theta_n$ in $a_n = 2\sqrt{p^n}\cos\theta_n$ is uniformly distributed in the interval $[-\pi, \pi]$ with respect to some measure as $n \to \infty$.

Proof: Let $L_n(T) = T^2 - a_n T + p^n$ be the characteristic polynomial of the Frobenious map, and consider $N_n = p^n + 1 - \alpha^n - \overline{\alpha}^n$, where $\alpha = (a_1 \pm \sqrt{a_1^2 - 4p})/2 = \sqrt{p} e^{i\theta_1}$ is a root of $L_1(T) = T^2 - a_1 T + p^n$. Then

$$N_n = p^n + 1 - \alpha^n - \overline{\alpha}^n = p^n + 1 - p^{n/2} e^{i\theta_1 n} - p^{n/2} e^{-i\theta_1 n} = p^n + 1 - 2p^{n/2} \cos(\theta_1 n)$$

Now by hypothesis $\theta = \theta_1 \ne r\pi$ is an irrational number. Hence, by Bohl's Lemma, (also known as Weil's Lemma) the sequence $\{ \theta n \mod 2\pi : n \geq 1 \}$ is uniformly distributed in the interval $[0, 2\pi]$ with respect to some measure. This in turn implies that the angle $\theta_n \equiv \theta n \mod 2\pi$ is uniformly distributed as $n \to \infty$. ∎

For rational $\theta = \theta_1 \in \mathbb{Q}$, or $\theta_1 = r\pi$, $r \in \mathbb{Q}$, it is quite clear that the subset $\{ \theta n \mod 2\pi : n \geq 1 \}$ is a finite subset of $[0, 2\pi]$.

Numerical data for $p = 2$ appear to show that the traces $a_n = 2\sqrt{p^n} \cos\theta_n$ are not uniformly distributed for some elliptic curves. For example, the elliptic curves $E_3 : y^2 + y = x^3 + x$ with $E_3(\mathbf{F}_2) \cong \mathbb{Z}_5$, and $E_4 : y^2 + y = x^3 + x + 1$ with $E_4(\mathbf{F}_2) \cong \mathbb{Z}_1$, the angles are $\theta_1 = \pi/2$ and $\theta_1 = -\pi/2$, respectively. In these exceptional cases, traces $a_n = 2\sqrt{2^n} \cos \pi n/2$, $n \geq 1$, are not uniformly distributed on the interval $2^n + 1 - 2\sqrt{2^n} \leq a_n \leq 2^n + 1 + 2\sqrt{2^n}$.

## 29. Frequency of the Values $a_p$

The cardinality of a subset of primes $S(t, x) = \#\{ p \leq x : p \nmid N_0 \text{ and } a_p = t \}$ is a topic of much interest in current research. The lower estimates for a few fixed values of $t$ are known. The frequency of supersingular elliptic curves is the same as the frequency of $a_p = 0$ which is exactly $1/2$ for CM curves. This turns out to be easier to establish.

**Theorem 67.** (Deuring 1941). A CM elliptic curve over the rational and conductor $N_0$ is supersingular if and only if the prime $p \nmid N_0$ is inert in the CM field of $E$ with finitely many exceptions of $p \nmid N_0$. Moreover, as $x$ tends to infinity,

$$\#\{ p \leq x : p \text{ does not divide } N_0 \text{ and } a_p = 0 \} = \frac{\pi(x)}{2} = \frac{x}{2\log(x)} + O\left(\frac{x}{\log^2(x)}\right). \tag{43}$$

On the other hand, Serre has demonstrated that the set of supersingular primes for nonCM elliptic curves has density zero, and the Lang-Trotter conjecture calls for

$$\#\{ p \leq x : E \text{ nonCM and supersingular at } p \} \sim cx^{1/2}/\log(x).$$

**Theorem 68.** (Deuring 1941). For any $t$ such that $0 \leq |t| \leq 2q^{1/2}$, there are $H(t^2 - 4q)$ elliptic curves over $\mathbf{F}_q$ of cardinality $\#E(\mathbf{F}_q) = q + 1 - t$ up to isomorphism.
The function $H(x)$ is the weighted class number, [SF87].

**Theorem 69.** ([FV96]). For any $\delta > 0$ there exists a real number $x_{E,\delta} > 0$ such that

$$\#\{ p \leq x : p \text{ does not divide } N \text{ and } a_p = t \} \geq \frac{\log\log\log(x)}{\log\log\log\log(x)^{1+\delta}}, \tag{44}$$

for $x > x_{E,\delta}$.

## 30. Zeta Functions of Curves of Genus $g \geq 1$

Let $\#C(\mathbf{F}_{q^n}) = q^n + 1 - a_n$ be the number of rational points on a curve $C$ of genus $g \geq 1$. The generating function of the integers $N_n = \#C(\mathbf{F}_{q^n})$ is called the *zeta function* of the curve. This power series is defined by

$$Z(T) = \exp\left(\sum_{n=1}^{\infty} \frac{N_n}{n} T^n\right) = \frac{L(T)}{(1-T)(1-qT)}. \tag{45}$$

The numerator is a polynomial $L(T) = c_{2g}T^{2g} + c_{2g-1}T^{2g-1} + \cdots + c_1 T + c_0$ of degree $2g$ with $g$ pairs of conjugate roots $\alpha_i \in \mathbb{C}$.

Some Properties of $L(T)$:
(i) $L(T) = (T - \alpha_1)(T - \alpha_2) \cdots (T - \alpha_{2g})$.
(ii) Conjugacy $\alpha_{i+1} = q/\alpha_i$, $1 \leq i \leq g$.
(iii) Uniform absolute value $|\alpha_i| = q^{1/2}$, $1 \leq i \leq 2g$.
(iv) $q$-Symmetric Coefficients: $c_{2g-i} = c_i q^{g-i}$, and $|a_i| \leq \binom{2g}{i} q^{i/2}$.
(v) The values $L(0) = q$, and $L(1) = \#C(\mathbf{F}_q)$.

The trace of $L(T)$ is determined using the formula,

$$a_n = \alpha_1^n + \overline{\alpha}_1^n + \alpha_2^n + \overline{\alpha}_2^n + \cdots + \alpha_g^n + \overline{\alpha}_g^n, \tag{46}$$

where the complex numbers $\alpha_1, \ldots, \alpha_g$ and their conjugates are the roots of $L(t)$ of absolute values $|\alpha_i| = q^{1/2}$, $1 \leq i \leq g$.

***Example 71.*** The first two characteristic polynomials are the following.
(1) Case of Elliptic Curves, $g = 1$. $L(T) = T^2 - t^n T + q^n$, where $t_n = \alpha_1^n + \overline{\alpha}_1^n$, and $\alpha_1^n \overline{\alpha}_1^n = q^n$.
The generating series is given by

$$Z(T) = \exp\left(\sum_{n=1}^{\infty} \frac{N_n}{n} T^n\right) = \frac{T^2 - t_n T + q^n}{(1-T)(1-qT)}. \tag{47}$$

(2) Case of Hyperelliptic Curves, $g = 2$. $L(T) = (T^2 + q) - aT(T^2 + q) + bT^2 = T^2 - (\alpha_1 + \alpha_2)T + q(T^2 - (\alpha_2 + \alpha_3)T + q$, where $a = \alpha_1^n + \overline{\alpha}_1^n + \alpha_2^n + \overline{\alpha}_2^n$, and $b = (\alpha_1^n + \overline{\alpha}_1^n)(\alpha_2^n + \overline{\alpha}_2^n)$, $|a| \leq q^{1/2}$, $b \leq 4q$. The generating series is given by

$$Z(T) = \exp\left(\sum_{n=1}^{\infty} \frac{N_n}{n} T^n\right) = \frac{T^4 - a_n T^3 + c_2 T^2 + c_1 T + q^{2n}}{(1-T)(1-qT)}. \tag{48}$$

The other way of writing it is $L(T) = T^4 - c_1 T^3 + c_2 T^2 - qc_1 T + q^2$, where $|c_1| \leq 4\sqrt{q}$, and $|c_2| \leq 6q$.

***Example 71.*** The zeta function of a supersingular elliptic curve is easier to determine since the roots of the characteristic polynomial are imaginary. Hence, $\alpha_1^n + \overline{\alpha}_1^n = 0$, and $\alpha_1^n \overline{\alpha}_1^n = q^n$ over any extension of $\mathbf{F}_q$ of odd degree $n$. The combined expression is

$$Z(T) = \frac{T^2 + q^n}{(1-T)(1-q^n T)}, \quad \text{and} \quad \#E(\mathbf{F}_{q^n}) = \begin{cases} q^n + 1 & \text{if } n \text{ is odd,} \\ q^n + 1 - 2(-q)^{n/2} & \text{if } n \text{ is even.} \end{cases} \quad (49)$$

**Lemma 72.** Let $\alpha \in \mathbb{Q}(\sqrt{d})$ be a complex integer, $d < 0$. Then the number $N_n = q^n + 1 - \alpha^n - \bar{\alpha}^n$ satisfies the norm equations

(1) $N_n = N(\alpha^n - 1)$,  (2) $N_n = N(\alpha^m - 1)N((\alpha^n - 1)/(\alpha^m - 1))$ for all $m \mid n$.

Each such integer $N_m = N(\alpha^m - 1)$ is a primitive divisor of the integer $N_n$.

## 31. L-Functions of Curves of Genus $g \geq 1$

The L-function of an elliptic curve is a global zeta function made up of the local zeta functions

$$Z_p(p^s) = 1 - a_p p^{-s} \text{ or } Z_p(p^s) = Z_p(p^s) = 1 - a_p p^{-s} + p^{1-2s}$$

at the primes $p$. The L-function is defined by

$$L(s) = \sum_{n=1}^{\infty} \frac{a_n}{n^s} = \prod_{\gcd(p, N_0) > 1} (1 - a_p p^{-s})^{-1} \prod_{\gcd(p, N) = 1} (1 - a_p p^{-s} + p^{1-2s})^{-1}, \quad (50)$$

where the first index runs over the prime divisors $p$ of $N_0$, and the second index runs over the primes $p$ relatively prime to $N_0$. The integer $N_0$ is the *conductor* of $E$.

**Definition 73.** The conductor $N_0$ of an elliptic curve $E : y^2 = f(x)$ is defined by

$$N_0 = \prod_{p \mid \Delta(E)} p^{f_p}, \quad \text{with} \quad f_p = \begin{cases} 0 & \text{if } E \text{ has good reduction at } p, \\ 1 & \text{if } E \text{ has multiplicative reduction at } p, \\ 2 & \text{if } E \text{ has additive reduction at } p \neq 2,3, \\ 2 + \delta_p & \text{if } E \text{ has additive reduction at } p = 2,3, \end{cases} \quad (51)$$

and $\delta_2 \leq 6$, $\delta_3 \leq 3$, see the literature.

This is the precise definition of the conductor of $E$, but for most practical purpose it is sufficient to take the discriminant $\Delta(E)$ or its squarefree part in place of $N_0$.

The twisted L-function is defined by $L(s, \chi) = \sum_{n \geq 1} a_n \chi(n) n^{-s}$, where $\chi(n) = \left(\frac{d}{n}\right)$ is the quadratic symbol.

The coefficients of the power series $f(z) = \sum_{n \geq 1} a_n q^n$, $(q = e^{i2\pi z})$, are completely determined by the generating set $\{a_p : p \text{ prime}\}$, and can be computed by means of the normalized relations:

(i) $a_1 = 1$,
(ii) $a_p a_q = a_{pq}$, $\gcd(p, q) = 1$,
(iii) $a_{p^n} = a_p a_{p^{n-1}} - p a_{p^{n-2}}$, $p$ prime.

***Example* 74.** The Fourier series of the elliptic curve $E: y^2 = x^3 + 1$. The coefficients were determined by either computing the traces $a_p$, $p$ prime, or using the normalized relations. For instance, $a_1 = 1$, $a_5 = 0$, $a_{5^2} = a_5 a_5 - 5a_1 = -5$, $a_{5^3} = a_5 a_{5^2} - 5a_5 = 0$, $a_{5^4} = a_5 a_{5^3} - 5a_{5^2} = 25$, etc. The power series is

$$f(z) = q - 4q^7 + 2q^{13} + 8q^{19} - 5q^{25} - 4q^{29} - 4q^{31} - 10q^{37} + 8q^{43} + 9q^{49} + \cdots .$$

**Note:** The coefficients $a_n \in \mathbb{Z}$ are the eigenvalues of the Hecke operator $T_n$ acting on the cusp forms of weight $k \geq 1$. Specifically, if $f(z) = \sum_{n=1}^{\infty} a_n q^n$ is the Fourier series expansion of a cusp form $f(z) \in S_k(\Gamma(N_0))$, then $T_n(f(z)) = a_n f(z)$. The weight 2 cusp form $f(z) \in S_2(\Gamma(N_0))$ of an elliptic curve with complex multiplication has infinitely many coefficients $a_p = 0$, (in fact, 50% of all the primes which are the primes of supersingular reduction). Moreover, there are arbitrary long sequences of consecutive vanishing values $a_n = a_{n+1} = \cdots = a_{n+m} = 0$.

***Theorem* 75.** ([AL05]). For every elliptic curve without complex multiplication, and any real numbers $x, y > 0$ such that $x^{51/134+\delta} < y$, $\delta > 0$, the following holds:

$\#\{x - y < n < x : a_n \neq 0\} \gg y$. And the maximal number of consecutive vanishing trace values $a_n = a_{n+1} = \cdots = a_{n+m} = 0$, satisfies $m \ll x^{51/134+\delta}$.

## 32. Average Size of the Groups of Points

The derivation of the mean value $\overline{\#E(\mathbf{F}_q)}$ of the size of group of points $\#E(\mathbf{F}_q)$ of an elliptic curve utilizes fairly advanced techniques. The analysis starts with the calculation of the (normalized) mean value $\overline{a_k(C)}$ of the trace of Frobenious $a_k(C) = \alpha_1^k + \overline{\alpha}_1^k + \cdots + \alpha_g^k + \overline{\alpha}_g^k$ of a curve $C$ of genus $g > 0$ over $\mathbf{F}_q$. The final expression is stated in terms of the Selberg trace formula

$$\sigma_k(T_p) = -1 - \frac{1}{2} \sum_{E/\mathbf{F}_q} \frac{\alpha^{k-1} - \overline{\alpha}^{k-1}}{\alpha - \overline{\alpha}}. \tag{52}$$

The summation index is over the set of elliptic curves $E$ over $\mathbf{F}_q$, and $T_p$ is the Hecke operator acting on cusp form of weight $k$ over the full group $SL(2, \mathbb{Z})$. The calculation uses the weights $1/\#\mathrm{Aut}(E)$ to remove duplicate curves.

***Theorem* 76.** ([BR01]). The (normalized) mean value of the trace $a_k(E)$ of an elliptic curve $E$ over $\mathbf{F}_q$ is given by

$$E\left(\frac{a_k}{q^{k/2}}\right) = \frac{1}{q} \sum_{E/\mathbf{F}_q} \frac{a_k(E)}{q^{k/2}} = \frac{\sigma_k(T_q) + 1}{q^{k/2}} - \frac{\sigma_{k+2}(T_q) + 1}{q^{k/2} + 1}. \tag{53}$$

The summation is taken over all the traces $a_k(E) = \alpha^k + \overline{\alpha}^k$ as the curve $E$ varies over $\mathbf{F}_q$.

***Corollary* 77.** ([BR01]). The mean value $\overline{\#E(\mathbf{F}_q)}$ of the size of group of points $\#E(\mathbf{F}_q)$ of an elliptic curve $E$ over $\mathbf{F}_q$ is $\overline{\#E(\mathbf{F}_q)} = q^k + O(q^{(k-1)/2+\varepsilon})$, $\varepsilon > 0$.

The last statement uses Deligne's estimate $O(q^{(k-1)/2+\varepsilon}) = \sigma_k(T_q)$ on the magnitude of the coefficients of the Fourier series of cusp forms of weight $k$. For small weight $k = 2, 4, 6, 8, 10, 12$, and $14$, the exact values are known:

$\sigma_0(T_q) = 0$    $\sigma_4(T_q) = 0$    $\sigma_8(T_q) = 0$    $\sigma_{12}(T_q) = \tau(q)$
$\sigma_2(T_q) = -q - 1$    $\sigma_6(T_q) = 0$    $\sigma_{10}(T_q) = 0$    $\sigma_{14}(T_q) = 0$

where $\tau(q)$ is the Ramanujan tau function. These values are used to compute the exact mean values of the first six even extensions:

$\overline{\#E(\mathbf{F}_q)} = q + o(q^{1/2})$,    $\overline{\#E(\mathbf{F}_{q^2})} = q^2 + q + 1 + 1/q$,

$\overline{\#E(\mathbf{F}_{q^k})} = q^k + 1/q$, $k = 4, 6, 8$,    $\overline{\#E(\mathbf{F}_{q^{10}})} = q^{10} + (\tau(q) + 1)/q$,

$\overline{\#E(\mathbf{F}_{q^{12}})} = q^{12} - \tau(q) + 1/q$.

The mean values of odd degree extensions appear to be unknown, for more details see [BR01].

**Note :** The Hecke operator $T_p$ and the Frobenius map $\sigma$ are linked by the relation $T_p = \sigma + p/\sigma$. The operator $T_p$ acts on the space $S_0(\Gamma_0(N_0))$ of cusp forms $f(z)$ of weight 2 for a curve of genus $g > 0$. This space is $g$-dimensional spanned by $f_1(z), f_2(z), \ldots, f_g(z)$. In the case of elliptic curves, $g = 1$, hence $S_0(\Gamma_0(N_0))$ has a single generator, see [HT03].

The coefficients $a_n \in \mathbb{Z}$ are the eigenvalues of the Hecke operator $T_n$ acting on the cusp form of weight 2. Specifically, if $f(z) \in S_2(\Gamma(N))$, then $T_n(f(z)) = a_n f(z)$.

## 33. Some Open Problems

(1) An algorithm to identify generators $g_1, \ldots, g_k$ of the group of points $E(\mathbf{F}_p) = <g_1> \times <g_2>$.

(2) A formula (analogous to Euler function $\varphi$) to enumerate the number of generators of the group of points $E(\mathbf{F}_p)$ of nonprime order. For groups of prime orders $N = \#E(\mathbf{F}_p)$, every point in $E(\mathbf{F}_p)$ has order $N$, so the number of generator is simply $N - 1$.

Acknowledgement: I would like to thank Professor Joseph H. Silverman for the comments, corrections and suggestions.